\documentclass[12pt,a4paper]{article}
\usepackage{multirow}
\usepackage{graphicx}
\usepackage{latexsym,amsmath,amsfonts,amscd,amssymb}
\usepackage{amsmath}
\usepackage{bm}
\usepackage{algorithm}
\usepackage{algorithmic}

\usepackage{epsfig}
\usepackage{epstopdf}
\usepackage{float}
\usepackage{changebar}
\usepackage{indentfirst}
\usepackage{verbatim}
\usepackage{color}
\usepackage{colordvi}
\allowdisplaybreaks[4]
\usepackage{appendix}
\usepackage{cite}



\newcommand{\eps}{\varepsilon}
\newtheorem{remark}{Remark}[section]

\begin{document}
\baselineskip=2pc
\begin{center}
{\Large \bf High order asymptotic preserving Hermite WENO fast sweeping method for the steady-state $S_{N}$ transport equation}
\end{center}
\centerline{Yupeng Ren\footnote{School of Mathematical Sciences, Xiamen University,
Xiamen, Fujian 361005, P.R. China. E-mail: ypren@stu.xmu.edu.cn. This work was carried out when Y. Ren was visiting Department of Mathematics, The Ohio State University under the support of the China Scholarship Council (CSC NO. 201906310077).},
Yulong Xing\footnote{Department of Mathematics, The Ohio State University, Columbus, OH 43210, USA.
 E-mail: xing.205@osu.edu. The work of Y. Xing is partially supported by the NSF grant DMS-1753581.},
 Dean Wang\footnote{Department of Mechanical and Aerospace Engineering, The Ohio State University, Columbus, OH 43210, USA.
 E-mail: wang.12239@osu.edu.} and Jianxian Qiu\footnote{School of Mathematical Sciences and Fujian Provincial
Key Laboratory of Mathematical Modeling and High-Performance
Scientific Computing, Xiamen University,
Xiamen, Fujian 361005, P.R. China. E-mail: jxqiu@xmu.edu.cn.
The work of J. Qiu is partially supported by NSFC grant 12071392.}}
\baselineskip=1.8pc
\vspace{1cm}
\centerline{\bf Abstract}
\bigskip
In this paper, we propose to combine the fifth order Hermite weighted essentially non-oscillatory (HWENO) scheme and fast sweeping method (FSM) for the solution of the steady-state $S_{N}$ transport equation in the finite volume framework. It is well-known that the $S_{N}$ transport equation asymptotically converges to a macroscopic diffusion equation in the limit of optically thick systems with small absorption and sources. Numerical methods which can preserve the asymptotic limit are referred to as asymptotic preserving methods. In the one-dimensional case, we provide the analysis to demonstrate the asymptotic preserving property of the high order finite volume HWENO method, by showing that its cell-edge and cell-average fluxes possess the thick diffusion limit.
Numerical results in both one- and two- dimensions are presented to validate its asymptotic preserving property.
A hybrid strategy to compute the nonlinear weights in the HWENO reconstruction is introduced to save computational cost.
Extensive one- and two-dimensional numerical experiments are performed to verify the accuracy, asymptotic
preserving property and positivity of the proposed HWENO FSM.

\vfill {\bf Key Words:} weighted essentially non-oscillatory (WENO) method; Hermite method; fast sweeping method; $S_{N}$ transport equation; asymptotic preserving property; diffusion limit


\pagenumbering{arabic}

\newpage

\section{Introduction}\setcounter{equation}{0}\setcounter{figure}{0}\setcounter{table}{0}
In this paper, we present a high order asymptotic preserving weighted essentially non-oscillatory (WENO) method for the steady-state transport equation, which can preserve the diffusion limit of the equation in the discrete setting. The radiative transport equation is a kinetic model which describes the scattering and absorbing of particles moving through a medium, and plays an important role in a wide range of scientific and engineering applications. The steady-state linear transport equation takes the form
\begin{subequations} \label{Eq:transport}
\begin{alignat}{2}
  &\Omega\cdot\nabla\psi(\textbf{x},\Omega)+\frac{\sigma_{t}}{\varepsilon}\psi(\textbf{x},\Omega)=
  \frac{1}{|\mathbb{S}^{d-1}|} \left( \left(\frac{\sigma_{t}}{\varepsilon}-\varepsilon\sigma_a\right){\phi}(\textbf{x})+\varepsilon Q(\textbf{x})\right),  \qquad
  	&& (\textbf{x},\Omega) \in D\times \mathbb{S}^{d-1}, \\
  &\psi(\textbf{x},\Omega)=G(\textbf{x},\Omega),  && (\textbf{x},\Omega) \in \Gamma^-,
\end{alignat}
\end{subequations}
where $D \subseteq \mathbb{R}^d$ (with $d=1,2,3$) is an open bounded domain.
When $d=3$, the set of propagation directions is the unit sphere $\mathbb{S}^{2}$ in $\mathbb{R}^3$,
When $d=1,2$, it becomes the projection of $\mathbb{S}^{2}$ onto $\mathbb{R}^d$, i.e., $\mathbb{S}^{1}$ is a unit disk if $d=2$ and
$\mathbb{S}^{0}$ is unit segment $[-1,1]$ if $d=1$.
$\Gamma^{-}=\left\{(\mathbf{x}, \mathbf{v}) \in \partial D \times \mathbb{S}^{d-1} ~ |~ \mathbf{n}(\mathbf{x}) \cdot \mathbf{v}<0\right\}$ is the incoming boundary, with $\mathbf{n}(\mathbf{x})$ being the unit outer normal vector at $\mathbf{x} \in \partial D$.
$\psi(\textbf{x},\Omega)$ denotes the angular intensity, and ${\phi}(\textbf{x})= \int_{\mathbb{S}^{d-1}} \psi d \Omega$ is the scalar flux representing the integral of $\psi$ over $\mathbb{S}^{d-1}$.
$\varepsilon$ is the scaling parameter, representing the ratio of a particle mean free path to a characteristic scale length of the system. 
$\sigma_{t}$, $\sigma_a$ and $\sigma_{s}$ are the non-dimensionalized total, absorption and scattering macroscopic cross section, which satisfy $\sigma_{a}=(\sigma_{t}-\sigma_{s})$. $Q(\textbf{x})$ is the external source function, and $G(\textbf{x},\Omega)$ is the given incoming flux on $\Gamma^-$.

It is well-known that when $\varepsilon$ is very small uniformly in the entire domain, the angular flux $\psi$ away from the boundary is nearly independent of the angular direction $\Omega$, and the transport model can be accurately approximated by a macroscopic diffusion equation that depends on the variable $\textbf{x}$ only \cite{1,2,3,4,5}. Asymptotic preserving (AP) numerical methods \cite{JinAP} refer to the methods that are accurate and robust in all regimes from transport dominated to diffusion dominated.
AP discretization of the transport equation \eqref{Eq:transport} reduces to a consistent and stable discretization of the macroscopic diffusion equation when $\varepsilon$ goes to zero.

There have been extensive studies on various AP numerical methods for solving the linear transport equation.
Larsen et al. first used asymptotic analysis to study the behavior of discrete transport solutions, and produced many important results on the relationship between the analytical and numerical solutions of the transport equation \cite{larsen1,larsen2,8,9,10}.  Larsen and others used the asymptotic expansion method to analyze
the behavior of several numerical schemes, such as the diamond difference
method \cite{L110,L11}, step difference method \cite{L110,L11}, the Lund-Wilson \cite{L112,L113} and Castor \cite{L114} methods in the
thick and intermediate regimes. Adams extended the asymptotic analysis to a complete family of discontinuous finite-element methods (DFEMs) and showed that some DFEM schemes do not possess the diffusion limit because the upwind numerical flux forces the scalar flux, and thus the angular flux, to be continuous across the mesh cells \cite{9}. Guermond and Kanschat proved by using functional analytic tools that a necessary and sufficient condition for the standard upwind discontinuous Galerkin approximation to converge to the correct limit solution in the diffusive regime is that the approximation space contains a linear space of continuous functions, and the restrictions of the functions of this space to each mesh cell contain the linear polynomials \cite{10}. Most recently, Wang has derived a theoretical result to determine the mesh size for a variety of finite difference scheme to achieve accurate results at the diffusion limit \cite{wangap}.
Finite volume Hermite WENO (HWENO) methods will be considered in this paper. WENO methods are a class of high order numerical methods for solving
the hyperbolic conservation laws, which yield very robust and non-oscillatory solutions near the shocks, and have been widely used in the applications.
Recently, high order HWENO methods, with a more compact reconstruction stencil, have also gained many attention in solving hyperbolic conservation laws.
The HWENO and WENO methods have the similar building block, and the major difference between them is that HWENO method uses both the unknown function and its first derivative (or first moment) in the reconstruction and update procedure. The HWENO scheme was first proposed as a robust limiter for the discontinuous Galerkin (DG) method in \cite{hwenolim1,hwenolim2}, thanks to the compact stencil required in its reconstruction step. In \cite{qiushuhwenohj2005,zhengdhj}, the HWENO scheme was extended to solve the Hamilton-Jacobi equation, and achieved very good numerical results. Compared with the standard WENO scheme, its boundary treatment is much simpler and the numerical error is observed to be smaller with the same mesh, as shown in \cite{qiushuhwenohj2005}. The HWENO scheme was later extended to solve the hyperbolic conservation laws in the finite difference \cite{lqfdhweno1,zhaohyhwenofd} and finite volume \cite{hwenolim1,hwenolim2} frameworks, and the same advantages have been observed.

In the past few decades, many efficient numerical solvers for the static hyperbolic conservation laws and Hamilton-Jacobi equations have been developed. Among them, one of the most popular method is the fast sweeping method (FSM) \cite{kaolf,qianzhangzhao,Tsai,Zhaofsm1}, which was first proposed by Bou\'{e} and Dupuis \cite{boue} to solve a deterministic control problem with quadratic running cost using Markov chain approximation. In \cite{Zhaofsm1}, a systematic way for solving the Eikonal equations using FSM was introduced by Zhao.
Later, many high order FSM have been developed to solving static Hamilton-Jacobi equations, in the framework of finite difference WENO \cite{zhang2006,huangnum,xiong} and finite element DG \cite{dgli,dgluo,dgwu,dgzhang} methods. In \cite{chenchou}, high order WENO FSM was proposed for solving the steady-state hyperbolic conservation laws with source terms. In
\cite{chenfix} and \cite{wuzhangfix}, FSM was combined with the fixed point iteration ideas to provide an efficient WENO solver for the steady-state hyperbolic conservation laws.

In this paper, we propose to combine the finite volume HWENO method with the fast sweeping technique, and apply them for the steady-state linear transport equation. In the angular discretization, we adopt the discrete ordinate ($S_N$) method, in which the angular variable is discretized into a finite number of directions, see \cite{larsen_SN} and the references therein for more discussion on $S_N$ method.
The main novel contribution of this paper is to present a class of high order AP methods, by demonstrating that the
proposed finite volume HWENO FSM preserves the asymptotic limit when $\varepsilon \rightarrow 0$.
Many high order AP methods have been studied for the linear transport equation in the literature, and most of their spatial discretizations are in the framework of DG method.
While DG methods enjoy many advantages including their robustness, flexibility and AP property (under certain conditions on the polynomial spaces, see \cite{9,10}), they are also known to be computationally expensive in multi dimensions when the polynomial degree becomes large.
In \cite{wangap}, it was shown that the original WENO method does not have AP property. We also investigated finite difference HWENO FSM and numerical results indicate that it is not AP.
Here we present a high order finite volume HWENO method (fifth-order HWENO is presented as an example, although the same idea can be extended to higher order if needed), which can be proven to have AP property following the similar approach in \cite{larsen2} to show the AP property of linear discontinuous (LD) method. The proposed method can also be viewed as the higher order extension of the LD method in one dimension, and that of the bilinear discontinuous finite element method \cite{ld2d} in multi dimensions.
In addition, we present a hybrid strategy to reduce the computational cost of evaluating the nonlinear weights in the HWENO reconstruction, which was shown to save about $50\%$ CPU time in the numerical tests.
In the two-dimensional case, we employ the dimension-by-dimension HWENO reconstruction procedure in the finite volume framework as in \cite{zhengdhj}, which can achieve the same essentially non-oscillatory property as the genuine two-dimensional strategy, and is easier to code than the latter one.
Both one- and two-dimensional algorithms have been studied, and extensive numerical examples are provided to confirm the AP property and robustness of the proposed methods.


The rest of the paper is organized as follows. In Section 2, we describe in detail the HWENO FSM for $S_{N}$ transport equation in one-dimensional (1D) case. The analysis of thick diffusion limit is also provided. In the Section 3, we introduce the HWENO FSM for multidimensional $S_{N}$ transport equation, and provide the flowchart of HWENO FSM in two-dimensional (2D) setting. The numerical examples are performed to demonstrate the high accuracy, positive and thick diffusion limit of our proposed schemes in Section 4. Some conclusion remarks are presented in Section 5.

\section{One-dimensional $S_{N}$ transport equation} \label{sec2}
\setcounter{equation}{0}\setcounter{figure}{0}\setcounter{table}{0}
In this section, we will present the HWENO FSM for 1D transport equation in the finite volume framework, and analyze the diffusion limit of the resulting method.

\subsection{Mathematical model and HWENO method}\label{sec:hweno1d}
The steady-state, monoenergetic, discrete ordinates $S_{N}$ transport equation in 1D
slab geometry $[0,L]$ with isotropic scattering takes the form
\begin{subequations}\label{1dsn}
\begin{align}
 &\mu_{m}\frac{d}{d x}\psi(x,\mu_{m})+\frac{\sigma_{t}}{\varepsilon}\psi(x,\mu_{m})=\frac{1}{2}\left(\frac{\sigma_t}{\varepsilon}-\varepsilon\sigma_a\right) \phi(x)+\frac{\varepsilon}{2}Q(x), \quad1\leq m\leq M,\\
 &\psi(0,\mu_{m})=f(\mu_{m}),\quad0<\mu_{m}\leq 1,\\
 &\psi(L,\mu_{m})=g(\mu_{m}),\quad-1\leq\mu_{m}<0,
\end{align}
\end{subequations}
where $\phi(x)=\int_{-1}^{1}\psi(x,\mu)d\mu=\sum\limits_{m=1}^{M}\psi(x,\mu_{m})\omega_{m}$ is the scalar flux, with $\omega_{m}$ being the Gaussian quadrature weights. $M$ is assumed to be an even integer in this paper, which means a symmetric quadrature set is used. The symmetric quadrature set $\{\mu_{m}, \omega_{m}\}$
satisfies
\begin{equation}\label{gaussweight}
\sum_{m=1}^{M}(\mu_{m})^{k}\omega_{m}=\begin{cases}
2,\quad &k=0,\\
0,\quad & \mathrm{for} ~k ~\mathrm{odd},\\
\frac{2}{k+1},\quad & \mathrm{for}~ k~ \mathrm{even},
\end{cases}
\end{equation}
where $k$ is an integer with $k\leq2M-1$.

\begin{figure}
\begin{center}
  \includegraphics[width=14cm]{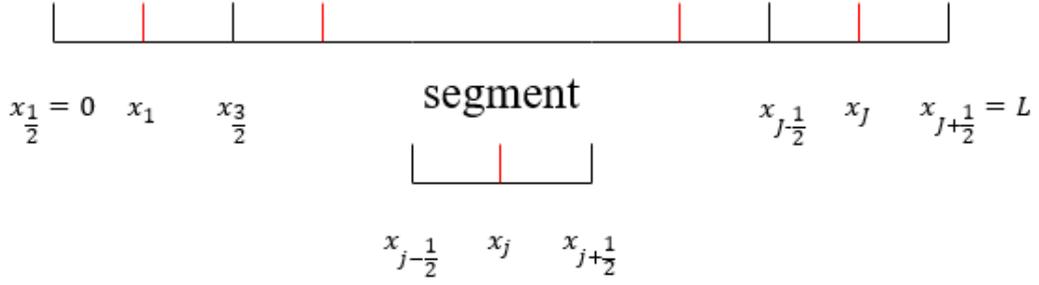}
\end{center}
\caption{The segment in one-dimensional case.}\label{segment}
\end{figure}

 We assume the computational domain $[0,L]$ has been divided into cells $I_{j}=[x_{j-\frac{1}{2}},x_{j+\frac{1}{2}}]$, $j=1,\cdots,J$ for all $m$. The cell center and the mesh size is denoted as $x_{j}=(x_{j-\frac{1}{2}}+x_{j+\frac{1}{2}})/2$ and $\Delta x_{j}=x_{j+\frac{1}{2}}-x_{j-\frac{1}{2}}$ as in Fig. \ref{segment}. Let $\psi_{m}(x)=\psi(x,\mu_m)$, we define
\begin{equation}
\psi_{m,j}=\frac{1}{\Delta x_{j}}\int_{I_{j}}\psi_{m}(x)dx,\quad
\widehat{\psi}_{m,j}=\frac{1}{\Delta x_{j}}\int_{I_{j}}\psi_{m}(x)\frac{x-x_{j}}{\Delta x_{j}}dx,
\end{equation}
as the numerical approximation to the average of angular flux and its first moment. Similarly, we can  define
\begin{align}
  &Q_{j}=\frac{1}{\Delta x_{j}}\int_{I_{j}}Q(x)dx,\quad
\widehat{Q}_{j}=\frac{1}{\Delta x_{j}}\int_{I_{j}}Q(x)\frac{x-x_{j}}{\Delta x_{j}}dx; \nonumber \\
  \displaystyle&\phi_{j}=\frac{1}{\Delta x_{j}}\int_{I_{j}}\phi(x)dx=\frac{1}{\Delta x_{j}}\int_{I_{j}}\int_{-1}^{1}\psi(x,\mu)d\mu dx=\sum\limits_{m=1}^{M}\psi_{m,j}\omega_{m};\label{phi_j}\\
  \displaystyle&\widehat{\phi}_{j}=\frac{1}{\Delta x_{j}}\int_{I_{j}}\phi(x)\frac{x-x_{j}}{\Delta x_{j}}dx=\frac{1}{\Delta x_{j}}\int_{I_{j}}\int_{-1}^{1}\psi(x,\mu)\frac{x-x_{j}}{\Delta x_{j}}d\mu dx=\sum\limits_{m=1}^{M}\widehat{\psi}_{m,j}\omega_{m}. \label{phi_jh}
\end{align}
Multiplying \eqref{1dsn} by $\frac{1}{\Delta x_{j}}$ and $\frac{x-x_{j}}{\Delta x_{j}^{2}}$ respectively, integrating on cell $I_{j}$ and applying integration by parts yield
\begin{subequations}
\begin{align}
&\frac{1}{\Delta x_{j}}\int_{I_{j}}(\mu_{m}\psi_{m}(x))_{x}dx+\frac{\sigma_{t,j}}{\varepsilon}\psi_{m,j}=\frac{1}{2}\left(\frac{\sigma_{t,j}}{\varepsilon}-\varepsilon\sigma_{a,j}\right)\sum\limits_{k=1}^{M}\psi_{k,j}\omega_{k}+\frac{\varepsilon}{2}Q_{j},\\
&\frac{1}{\Delta x_{j}}\left(\frac{1}{2}\left(\mu_{m}\psi_{m,j+\frac{1}{2}}+\mu_{m}\psi_{m,j-\frac{1}{2}}\right)-\mu_{m}\psi_{m,j}\right)+\sigma_{t,j}\widehat{\psi}_{m,j}=\frac{1}{2}\left(\frac{\sigma_{t,j}}{\varepsilon}-\varepsilon\sigma_{a,j}\right)\sum\limits_{k=1}^{M}\widehat{\psi}_{k,j}\omega_{k}+\frac{\varepsilon}{2}\widehat{Q}_{j},\\
&\psi(0,\mu_{m})=f(\mu_{m}),\quad0<\mu_{m}\leq 1,\\
&\psi(L,\mu_{m})=g(\mu_{m}),\quad-1\leq\mu_{m}<0,
\end{align}
\end{subequations}
where the following equality with $F_m(\psi)=\mu_{m}\psi_{m}(x)$
\begin{equation}
    \int_{I_{j}}(F_m(\psi))_{x}\frac{x-x_{j}}{\Delta x_{j}}dx=\frac{1}{2}\left(F_{m}(x_{j+\frac{1}{2}}^-)+F_{m}(x_{j-\frac{1}{2}}^+)\right)-\mu_{m}\psi_{m,j}
\end{equation}

is used in the derivation of the second equation.
The HWENO numerical discretizations of equation \eqref{1dsn} are now given by
\begin{equation}\label{discre1d}
\begin{aligned}
& \frac{1}{\Delta x_{j}}\left(\widehat{F}_{m,j+\frac{1}{2}}-\widehat{F}_{m,j-\frac{1}{2}}\right)+\frac{\sigma_{t,j}}{\varepsilon}\psi_{m,j}=\frac{1}{2}\left(\frac{\sigma_{t,j}}{\varepsilon}-\varepsilon\sigma_{a,j}\right)\sum\limits_{k=1}^{M}\psi_{k,j}\omega_{k}+\frac{\varepsilon}{2}Q_{j},\\
& \frac{1}{\Delta x_{j}}\left(\frac{1}{2}(\widehat{F}_{m,j+\frac{1}{2}}+\widehat{F}_{m,j-\frac{1}{2}})-\mu_{m}\psi_{m,j}\right)+\frac{\sigma_{t,j}}{\varepsilon}\widehat{\psi}_{m,j}=\frac{1}{2}\left(\frac{\sigma_{t,j}}{\varepsilon}-\varepsilon\sigma_{a,j}\right)\sum\limits_{k=1}^{M}\widehat{\psi}_{k,j}\omega_{k}+\frac{\varepsilon}{2}\widehat{Q}_{j},
\end{aligned}
\end{equation}
with $\widehat{F}_{m,j\pm\frac{1}{2}}$ being the numerical fluxes to be specified, and
$\psi_{m,\frac{1}{2}}=f_{m}$ if $\mu_{m}>0$ and $\psi_{m,J+\frac{1}{2}}=g_{m}$ if $\mu_{m}<0$.
In this paper, we use the Godunov numerical flux, taking the form
\begin{equation}
\widehat{F}_{m,j+\frac{1}{2}}=
\begin{cases}
\displaystyle \min\limits_{\psi_{m,j+\frac{1}{2}}^{-}\leq \psi\leq \psi_{m,j+\frac{1}{2}}^{+}}\mu_{m}\psi_{m}, \quad \text{if}~ \psi_{m,j+\frac{1}{2}}^{-}\leq \psi_{m,j+\frac{1}{2}}^{+};\\
\displaystyle  \max\limits_{\psi_{m,j+\frac{1}{2}}^{+}\leq \psi\leq \psi_{m,j+\frac{1}{2}}^{-}}\mu_{m}\psi_{m}, \quad \text{if} ~\psi_{m,j+\frac{1}{2}}^{-}> \psi_{m,j+\frac{1}{2}}^{+}.
\end{cases}
\end{equation}
The Godonov flux can be further simplified for our linear flux, which leads to the following two cases:
\begin{itemize}
  \item If $\mu_{m}>0$, the HWENO scheme \eqref{discre1d} becomes
  \begin{equation}\label{discre1d1}
\begin{aligned}
&  \frac{\mu_{m}}{\Delta x_{j}}\left(\psi_{m,j+\frac{1}{2}}^{-}-\psi_{m,j-\frac{1}{2}}^{-}\right)+\frac{\sigma_{t,j}}{\varepsilon}\psi_{m,j}=\frac{1}{2}\left(\frac{\sigma_{t,j}}{\varepsilon}-\varepsilon\sigma_{a,j}\right)\sum\limits_{k=1}^{M}\psi_{k,j}\omega_{k}+\frac{\varepsilon}{2}Q_{j},\\
&  \frac{\mu_{m}}{\Delta x_{j}}\left(\frac{1}{2}(\psi_{m,j+\frac{1}{2}}^{-}+\psi_{m,j-\frac{1}{2}}^{-})-\psi_{m,j}\right)+\frac{\sigma_{t,j}}{\varepsilon}\widehat{\psi}_{m,j}=\frac{1}{2}\left(\frac{\sigma_{t,j}}{\varepsilon}-\varepsilon\sigma_{a,j}\right)\sum\limits_{k=1}^{M}\widehat{\psi}_{k,j}\omega_{k}+\frac{\varepsilon}{2}\widehat{Q}_{j},
\end{aligned}
\end{equation}
  with $\psi_{m,\frac{1}{2}}=f_{m}$.

  \item If $\mu_{m}<0$, the HWENO scheme \eqref{discre1d} becomes
    \begin{equation}\label{discre1d2}
\begin{aligned}
&  \frac{\mu_{m}}{\Delta x_{j}}\left(\psi_{m,j+\frac{1}{2}}^{+}-\psi_{m,j-\frac{1}{2}}^{+}\right)+\frac{\sigma_{t,j}}{\varepsilon}\psi_{m,j}=\frac{1}{2}\left(\frac{\sigma_{t,j}}{\varepsilon}-\varepsilon\sigma_{a,j}\right)\sum\limits_{k=1}^{M}\psi_{k,j}\omega_{k}+\frac{\varepsilon}{2}Q_{j},\\
&  \frac{\mu_{m}}{\Delta x_{j}}\left(\frac{1}{2}(\psi_{m,j+\frac{1}{2}}^{+}+\psi_{m,j-\frac{1}{2}}^{+})-\psi_{m,j}\right)+\frac{\sigma_{t,j}}{\varepsilon}\widehat{\psi}_{m,j}=\frac{1}{2}\left(\frac{\sigma_{t,j}}{\varepsilon}-\varepsilon\sigma_{a,j}\right)\sum\limits_{k=1}^{M}\widehat{\psi}_{k,j}\omega_{k}+\frac{\varepsilon}{2}\widehat{Q}_{j},
\end{aligned}
\end{equation}
   with $\psi_{m,J+\frac{1}{2}}=g_{m}$.
\end{itemize}

Next, we will present the HWENO reconstruction procedure to evaluate the high order interface value approximations $\psi_{m,j\pm\frac{1}{2}}^{\mp}$ from the cell average values $\psi_{m,j}$ and the first order moments $\widehat{\psi}_{m,j}$. For ease of presentation,
we assume the mesh is uniform, i.e. $\Delta x_{j}=\Delta x$ for all $j$. The detailed procedure of the HWENO reconstruction is summarized as follows:

 1. Based on three small stencils $S_{0}=\{I_{j-1},I_{j}\}$, $S_{1}=\{I_{j},I_{j+1}\}$, $S_{2}=\{I_{j-1},I_{j},I_{j+1}\}$, and a bigger stencil $T=\{S_{0},S_{1},S_{2}\}$, we construct three Hermite cubic polynomials $p_{0}(x),~p_{1}(x),~p_{2}(x)$, and a fifth-order polynomial $q(x)$ such that
\begin{align*}
\frac{1}{\Delta x}\int_{I_{j+i}}p_{0}(x)dx=\psi_{m,j+i}, \quad &\frac{1}{\Delta x}\int_{I_{j+i}}p_{0}(x)\frac{x-x_{j+i}}{\Delta x}dx=\widehat{\psi}_{m,j+i},  & i=-1,0,\\
\frac{1}{\Delta x}\int_{I_{j+i}}p_{1}(x)dx=\psi_{m,j+i},\quad &\frac{1}{\Delta x}\int_{I_{j+i}}p_{1}(x)\frac{x-x_{j+i}}{\Delta x}dx=\widehat{\psi}_{m,j+i},  & i=0,1,\\
\frac{1}{\Delta x}\int_{I_{j+i}}p_{2}(x)dx=\psi_{m,j+i}, \quad&\frac{1}{\Delta x}\int_{I_{j}}p_{2}(x)\frac{x-x_{j}}{\Delta x}dx=\widehat{\psi}_{m,j},  & i=-1,0,1,\\
\frac{1}{\Delta x}\int_{I_{j+i}}q(x)dx=\psi_{m,j+i}, \quad&\frac{1}{\Delta x}\int_{I_{j+i}}q(x)\frac{x-x_{j+i}}{\Delta x}dx=\widehat{\psi}_{m,j+i},  & i=-1,0,1.
\end{align*}
Only the values of these polynomials at the cell interfaces $x=x_{j\pm\frac{1}{2}}$ are needed, and they take the form
\begin{subequations}
\begin{align}
\displaystyle
&p_{0}(x_{j-\frac{1}{2}}^{+})=\frac{1}{2}\psi_{m,j-1}+\frac{1}{2}\psi_{m,j}+2\widehat{\psi}_{m,j-1}-2\widehat{\psi}_{m,j}; \label{linear11} \\
\displaystyle&p_{0}(x_{j+\frac{1}{2}}^{-})=\frac{3}{4}\psi_{m,j-1}+\frac{1}{4}\psi_{m,j}+\frac{7}{2}\widehat{\psi}_{m,j-1}+\frac{23}{2}\widehat{\psi}_{m,j}; \label{linear12} \\
\displaystyle
&p_{1}(x_{j-\frac{1}{2}}^{+})=\frac{1}{4}\psi_{m,j}+\frac{3}{4}\psi_{m,j+1}-\frac{23}{2}\widehat{\psi}_{m,j}-\frac{7}{2}\widehat{\psi}_{m,j+1}; \label{linear13} \\
\displaystyle&p_{1}(x_{j+\frac{1}{2}}^{-})=\frac{1}{2}\psi_{m,j}+\frac{1}{2}\psi_{m,j+1}+2\widehat{\psi}_{m,j}-2\widehat{\psi}_{m,j+1}; \label{linear14} \\
\displaystyle
&p_{2}(x_{j-\frac{1}{2}}^{+})=\frac{7}{66}\psi_{m,j-1}+\frac{5}{6}\psi_{m,j}+\frac{2}{33}\psi_{m,j+1}-\frac{60}{11}\widehat{\psi}_{m,j}; \label{linear15} \\
\displaystyle&p_{2}(x_{j+\frac{1}{2}}^{-})=\frac{2}{33}\psi_{m,j-1}+\frac{5}{6}\psi_{m,j}+\frac{7}{66}\psi_{m,j+1}+\frac{60}{11}\widehat{\psi}_{m,j}; \label{linear16} \\
\displaystyle
&q(x_{j-\frac{1}{2}}^{+})=\frac{8}{27}\psi_{m,j-1}+\frac{7}{12}\psi_{m,j}+\frac{13}{108}\psi_{m,j+1}+\frac{28}{27}\widehat{\psi}_{m,j-1}-\frac{241}{54}\widehat{\psi}_{m,j}-\frac{25}{54}\widehat{\psi}_{m,j+1}; \label{linear1}\\
\displaystyle
&q(x_{j+\frac{1}{2}}^{-})=\frac{13}{108}\psi_{m,j-1}+\frac{7}{12}\psi_{m,j}+\frac{8}{27}\psi_{m,j+1}+\frac{25}{54}\widehat{\psi}_{m,j-1}+\frac{241}{54}\widehat{\psi}_{m,j}-\frac{28}{27}\widehat{\psi}_{m,j+1};\label{linear2}
\end{align}
\end{subequations}

2. For each small stencil $S_{k}, k=0,1,2$, we compute the smooth indicators respectively, which measure the smoothness of the reconstructed polynomials $p_{k}(x),~k=0,1,2$, in the target cell $I_{j}$. The smaller the indicator is, the smoother the polynomial is in the target cells. Generally speaking, the smooth indicators are defined as
\begin{equation}\label{indicator}
  \beta_{k}=\sum\limits_{l=1}^{3}\int_{I_{j}}\Delta x^{2l-1}\left(\frac{\partial ^{l}}{\partial x^{l}}p_{k}(x)\right)^{2},	\qquad k=0,1,2,
\end{equation}
and their specific expressions are given by
\begin{align}
\displaystyle
&\beta_{0}=\frac{1}{16}(\psi_{m,j}-\psi_{m,j-1}-54\widehat{\psi}_{m,j}-6\widehat{\psi}_{m,j-1})^2+\frac{39}{16}(-5\psi_{m,j-1}+5\psi_{m,j}-38\widehat{\psi}_{m,j}-22\widehat{\psi}_{m,j-1})^2\nonumber \\
\displaystyle&\qquad+\frac{3905}{16}(-\psi_{m,j-1}+\psi_{m,j}-6\widehat{\psi}_{m,j}-6\widehat{\psi}_{m,j-1})^2, \nonumber\\
\displaystyle&\beta_{1}=\frac{1}{16}(\psi_{m,j}-\psi_{m,j+1}+54\widehat{\psi}_{m,j}+6\widehat{\psi}_{m,j+1})^2+\frac{39}{16}(-5\psi_{m,j+1}+5\psi_{m,j}+38\widehat{\psi}_{m,j}+22\widehat{\psi}_{m,j+1})^2 \nonumber\\
\displaystyle&\qquad+\frac{3905}{16}(-\psi_{m,j+1}+\psi_{m,j}+6\widehat{\psi}_{m,j}+6\widehat{\psi}_{m,j+1})^2, \nonumber\\
\displaystyle&\beta_{2}=\frac{1}{484}(-\psi_{m,j-1}+\psi_{m,j+1}+240\widehat{\psi}_{m,j})^2+\frac{13}{12}(-\psi_{m,j-1}+2\psi_{m,j}-\psi_{m,j+1})^2 \notag \\
\displaystyle&\qquad+\frac{355}{44}(-\psi_{m,j+1}+\psi_{m,j-1}+24\widehat{\psi}_{m,j})^2.\nonumber
\end{align}
In this paper, we follow the approach in \cite{wanglfweno}, and use $\beta'_{m}=\tau_{m}\beta_{m}$ as the smoothness indicator, where
\begin{equation*}
  \begin{split}
     \tau_{0}=&\max\left[|\sigma_{t,j+1}-\sigma_{t,j}|,|\sigma_{s,j+1}-\sigma_{s,j}|\right]\Delta x, \\
       \tau_{1}=&\max\left[|\sigma_{t,j}-\sigma_{t,j-1}|,|\sigma_{s,j}-\sigma_{s,j-1}|\right]\Delta x
   \end{split}
\end{equation*}
and $\tau_{2}=\max[\tau_{0},\tau_{1}]$. These parameters are introduced in \cite{wanglfweno} to estimate the local material heterogeneity, and for the steady-state linear problem studied in this paper, it is known that the discontinuity will appear only at the location when the heterogeneity occurs.

\indent 3. We compute the linear weights, denoted by $\gamma_k(x_{j\pm\frac{1}{2}}^\mp),~k=0,1,2$, satisfying
\begin{equation*}
  q(x_{j\pm\frac{1}{2}}^{\mp})=\sum\limits_{k=0}^{2}\gamma_{k}(x_{j\pm\frac{1}{2}}^{\mp})p_{k}(x_{j\pm\frac{1}{2}}^{\mp})
\end{equation*}
in the smooth region, which leads to the values
\begin{equation*}
  \begin{split}
     \gamma_{0}(x_{j-\frac{1}{2}}^{+})=\frac{14}{27},\quad \gamma_{1}(x_{j-\frac{1}{2}}^{+})=\frac{25}{189},\quad\gamma_{2}(x_{j-\frac{1}{2}}^{+})=\frac{22}{63}, \\
     \gamma_{0}(x_{j+\frac{1}{2}}^{-})=\frac{25}{189},\quad \gamma_{1}(x_{j+\frac{1}{2}}^{-})=\frac{14}{27},\quad\gamma_{2}(x_{j+\frac{1}{2}}^{-})=\frac{22}{63}.
   \end{split}
\end{equation*}

4. Combining the linear weights, smoothness indicators, one can evaluate the nonlinear weights by
\begin{equation}\label{nonweight}
\overline{\omega}_{k}(x_{j\pm\frac{1}{2}}^{\mp}) =\frac{\widetilde{\omega}_{k}(x_{j\pm\frac{1}{2}}^{\mp})}{\sum\limits_{\ell}\widetilde{\omega}_{\ell}(x_{j\pm\frac{1}{2}}^{\mp})},\qquad
\widetilde{\omega}_{k}(x_{j\pm\frac{1}{2}}^{\mp})=\frac{\gamma_{k}(x_{j\pm\frac{1}{2}}^{\mp})}{(\beta_{k}'+\widetilde{\varepsilon})^{2}},\qquad k=0,1,2,
\end{equation}
where $\widetilde{\varepsilon}$ is a small positive number to avoid the denominator becoming zero, and is taken as $\widetilde{\varepsilon}=10^{-6}$. The actual HWENO approximations of the cell interface values take the form
\begin{equation}\label{nonlinear}
  \psi_{_{m,j\pm\frac{1}{2}}}^{\mp} = \sum_{k=0}^{2}\overline{\omega}_{k}(x_{j\pm\frac{1}{2}}^{\mp})p_{k}(x_{j\pm\frac{1}{2}}^{\mp}),
\end{equation}
with $p_k(x_{j\pm\frac{1}{2}}^{\mp})$ defined in \eqref{linear11}-\eqref{linear16}.

\begin{remark} \label{remark2.1}
It is not difficult to observe from \eqref{indicator} that, if both $\sigma_{t}$ and $\sigma_{s}$ are constants in the big stencil $T$, the corresponding $\beta_{m}'$ equals to zero. Therefore, we can replace the nonlinear HWENO reconstruction \eqref{nonlinear} by the following linear approximation in \eqref{linear1}-\eqref{linear2}
\begin{subequations}
\begin{align}
&\psi_{m,j-\frac{1}{2}}^{+}=\frac{8}{27}\psi_{m,j-1}+\frac{7}{12}\psi_{m,j}+\frac{13}{108}\psi_{m,j+1}+\frac{28}{27}\widehat{\psi}_{m,j-1}-\frac{241}{54}\widehat{\psi}_{m,j}-\frac{25}{54}\widehat{\psi}_{m,j+1}; \label{linear21}\\
&\psi_{m,j+\frac{1}{2}}^{-}=\frac{13}{108}\psi_{m,j-1}+\frac{7}{12}\psi_{m,j}+\frac{8}{27}\psi_{m,j+1}+\frac{25}{54}\widehat{\psi}_{m,j-1}+\frac{241}{54}\widehat{\psi}_{m,j}-\frac{28}{27}\widehat{\psi}_{m,j+1}.\label{linear22}
\end{align}
\end{subequations}
This hybrid strategy is valid, since the shock will not appear for this steady-state linear equation
in the region when $\sigma_t$ and $\sigma_s$ are both constants. From the numerical results in the Section \ref{sec-num}, it can be observed that this strategy can save about $50\%$ CPU time. 
\end{remark}

\subsection{Fast sweeping idea to solve the global linear system}

The proposed HWENO scheme for the linear transport equation takes the form of \eqref{discre1d1} or \eqref{discre1d2}, combined with the HWENO reconstruction of $\psi_{m,j\pm\frac{1}{2}}^{\mp}$. This is a large system involving the flux term (coupling in $x$ direction) on the left side and the summation term (coupling in $\Omega$ direction) on the right side. The fast sweeping idea is adopted to solve this system efficiently.
Let us first denote the right-hand side term of the two equations in \eqref{discre1d1} as $S_j$ and $\widehat{S}_j$, respectively.
We summarize the flowchart of HWENO FSM for $S_{N}$ equations in 1D as follows
and refer to \cite{zhang2006,Zhaofsm1} for details of the FSM.

\noindent \textbf{Step 1}. \emph{Initialization}:
 We take $0$ as the initial guess of the unknowns $\psi_{m,j}$ and $\widehat{\psi}_{m,j}$ for all $m$ and $j$, and evaluate $S_j$ and $\widehat{S}_j$.

\noindent\textbf{Step 2}. \emph{Gauss-Seidel iteration with alternating sweep}. We sweep the whole domain with the following two alternating orderings repeatedly for each $m$: 

(I) $j=1\rightarrow J$: if $\mu_{m}>0$, solve the system \eqref{discre1d1} for each $j$ from left to right. After updating the approximation $\psi_{m,j}$ and $\widehat{\psi}_{m,j}$ in the cell $I_{j}$, we can apply HWENO reconstruction to obtain the cell-edge flux $\psi_{m,j+\frac{1}{2}}^{-}$ based on the most updated values of $\psi_{m,j+i}$ and $\widehat{\psi}_{m,j+i}$ ($i=-1,0,1$).

(II) $j=J\rightarrow1$: if $\mu_{m}<0$, solve the system \eqref{discre1d2} for each $j$ from right to left. After updating the approximation $\psi_{m,j}$ and $\widehat{\psi}_{m,j}$ in the cell $I_{j}$, we can apply HWENO reconstruction to obtain the cell-edge flux $\psi_{m,j-\frac{1}{2}}^{+}$ based on the most updated values of $\psi_{m,j+i}$ and $\widehat{\psi}_{m,j+i}$ ($i=-1,0,1$).

At the boundary of the computational domain, high order extrapolations are used to compute the values at the ghost cells, which are needed for the HWENO reconstruction near the boundary. After repeating this process for all $m$ directions, we can compute scalar flux $\phi_j$, $\widehat{\phi}_j$ from $\psi_{m,j}$, $\widehat{\psi}_{m,j}$ via Gauss quadrature \eqref{phi_j}-\eqref{phi_jh}, and update $S_j$, $\widehat{S}_j$. This completes one Gauss-Seidel iteration.

\noindent \textbf{Step 3}. \emph{Convergence}: Repeat the Gauss-Seidel iteration until the convergence criteria is satisfied. In this paper, if the scalar flux satisfies
$$\delta=||\phi^{new}-\phi^{old}||_{L_{1}}<10^{-14},$$
for two consecutive iteration steps, we stop the iteration.

The pseudo code of \textbf{Step 2} is presented in Algorithm \ref{Algorithm1}, where the superscript $n$ indicates the results in ``$n$-th" iteration.

\begin{algorithm}[h]
\caption{The Gauss-Seidel iteration of evaluating the scalar flux $\psi^{n+1}$ from $\psi^{n}$}
\label{Algorithm1}
\begin{algorithmic}
\REQUIRE  The values of $\psi_{m,j}, S_j$ and $\widehat{S}_j$ after $n$-th iteration.
\FOR{$m=1$ \TO M}
{
\IF{$\mu_{m}>0$}
\STATE $\psi_{m,\frac{1}{2}} \gets f_{m}$
{
\FOR {$j=1$ \TO J}
\STATE Solve \eqref{discre1d1} to obtain
$\psi_{m,j}$ and $\widehat{\psi}_{m,j}$\\
\STATE Compute cell-edge $\psi_{m,j+\frac{1}{2}}=\psi_{m,j+\frac{1}{2}}^{-}$ by HWENO reconstruction
\ENDFOR
\STATE High order extrapolation are used to compute the values at the ghost cells
}
\ELSE
\STATE $\psi_{m,J+\frac{1}{2}} \gets g_{m}$
\FOR {$j=J$ \TO 1}
\STATE Solve \eqref{discre1d2} to obtain $\psi_{m,j}$ and $\widehat{\psi}_{m,j}$
\STATE Compute cell-edge $\psi_{m,j-\frac{1}{2}}=\psi_{m,j-\frac{1}{2}}^{+}$ by HWENO reconstruction
\ENDFOR
\STATE High order extrapolation are used to compute the values at the ghost cells
\ENDIF
}
\ENDFOR
\FOR{$j=1$ \TO J}
\STATE $\phi_{j}=\sum\limits_{m=1}^{M}\psi_{m,j}\omega_{m}$, \quad $\widehat{\phi}_{j}=\sum\limits_{m=1}^{M}\widehat{\psi}_{m,j}\omega_{m}$
\STATE $S_{j}=\frac{1}{2}\left(\frac{\sigma_{t,j}}{\varepsilon}-\varepsilon\sigma_{a,j}\right)\phi_{j}+\frac{\varepsilon}{2}Q_{j}$, \quad
$\widehat{S}_{j}=\frac{1}{2}\left(\frac{\sigma_{t,j}}{\varepsilon}-\varepsilon\sigma_{a,j}\right)\widehat{\phi}_{j}+\frac{\varepsilon}{2}\widehat{Q}_{j}$
\ENDFOR
\IF{$\delta=||\phi^{n+1}-\phi^{n}||_{L_{1}}<10^{-14}$}
\STATE Stop iterate
\ELSE
\STATE Return to the top and continue the iteration
\ENDIF
\end{algorithmic}
\end{algorithm}

\begin{remark}\label{BC}
The high order extrapolations are used to evaluate the values at the ghost cells.
Because the stencil of HWENO method is more compact than that of WENO method, we need only one ghost cell in left boundary and right boundary of computational domain, respectively.
\end{remark}

\subsection{Thick diffusion limit}

One focus of the proposed HWENO method is its AP property when $\varepsilon$ is small.
In this subsection, we will provide the mathematical analysis to study the thick diffusion limit of the HWENO method. It will be showed that the cell-edge and cell-average fluxes possess the thick diffusion limit, and the HWENO method is very accurate for problems with anisotropic boundary fluxes. The detailed analysis is inspired by that of the LD method in \cite{larsen2}.

%
It is known that when $\varepsilon\rightarrow0$, the solution of the 1D linear transport equation \eqref{1dsn} satisfies \cite{larsen1,larsen2}
\begin{equation*}
  \psi(x,\mu)=\frac{\phi(x)}{2}+O(\varepsilon),
\end{equation*}
where $\phi(x)$ is the solution of the diffusion equation
\begin{equation}\label{diffusion}
  -\frac{d}{dx}\frac{1}{3\sigma_{t}}\frac{d}{dx}\phi+\sigma_{a}\phi=Q.
\end{equation}
with appropriate boundary conditions, and we refer to \cite{larsen2} for more discussions on this.

We will analyze the asymptotic diffusion limit of the HWENO FSM with linear reconstruction, and verify that it is a consistent approximation of the diffusion equation \eqref{diffusion}. Let us rewrite the HWENO FSM \eqref{discre1d1}-\eqref{discre1d2} in the following asymptotic form
\begin{equation}\label{1dHWENO_AP}
\begin{aligned}
& \frac{\mu_{m}}{\Delta x_{j}}\left(\psi_{m,j+\frac{1}{2}}-\psi_{m,j-\frac{1}{2}}\right)+\frac{\sigma_{t,j}}{\varepsilon}\psi_{m,j}=\frac{1}{2}\left(\frac{\sigma_{t,j}}{\varepsilon}-\varepsilon\sigma_{a,j}\right)\sum\limits_{k=1}^{M}\psi_{k}\omega_{k}+\frac{\varepsilon Q_{j}}{2},\\
& \frac{\mu_{m}}{\Delta x_{j}}\left(\frac{1}{2}(\psi_{m,j+\frac{1}{2}}+\psi_{m,j-\frac{1}{2}})-\psi_{m,j}\right)+\frac{\sigma_{t,j}}{\varepsilon}\widehat{\psi}_{m,j}=\frac{1}{2}\left(\frac{\sigma_{t,j}}{\varepsilon}-\varepsilon\sigma_{a,j}\right)\sum\limits_{k=1}^{M}\widehat{\psi}_{k}\omega_{k}+\frac{\varepsilon \widehat{Q}_{j}}{2},\\
\end{aligned}
\end{equation}
where we ignore the ``$\pm$" sign in the numerical fluxes for simplicity, and $\psi_{m,\frac{1}{2}}=f_{m}$ if $\mu_{m}>0$, $\psi_{m,J+\frac{1}{2}}=g_{m}$ if $\mu_{m}<0$. To perform the asymptotic analysis, we start by introducing the following ansatz
\begin{equation*}
\psi_{m} = \sum\limits_{k=0}^{\infty}\varepsilon^{k}\psi_{m}^{(k)}
\end{equation*}
for both the cell-edge fluxes $\psi_{m,j+\frac{1}{2}}$ and cell-average fluxes $\psi_{m,j}$ in \eqref{discre1d}. After plugging this ansatz into the HWENO method
\eqref{1dHWENO_AP}, we collect the equations with different orders of $\varepsilon$ and present the detailed proof of the AP property in four steps summarized as follows.

\underline{\textbf{Step 1, $O(\varepsilon^{-1})$ equations:}} The two $O(\varepsilon^{-1})$ equations are
\begin{equation}
\begin{split}
&\sigma_{t,j}\left(\psi_{m,j}^{(0)}-\frac{1}{2}\sum\limits_{k=1}^{M}\psi_{k,j}^{(0)}\omega_{k}\right)=0,\\
&\sigma_{t,j}\left(\widehat{\psi}_{m,j}^{(0)}-\frac{1}{2}\sum\limits_{k=1}^{M}\widehat{\psi}_{k,j}^{(0)}\omega_{k}\right)=0.
\end{split}
\end{equation}
These equations have isotropic solutions
\begin{equation}\label{eps0}
\begin{split}
\psi_{m,j}^{(0)}=\frac{1}{2}\sum\limits_{k=1}^{M}\psi_{k,j}^{(0)}\omega_{k}=\frac{1}{2}\phi_{j}^{(0)},\\
\widehat{\psi}_{m,j}^{(0)}=\frac{1}{2}\sum\limits_{k=1}^{M}\widehat{\psi}_{k,j}^{(0)}\omega_{k}=\frac{1}{2}\widehat{\phi}_{j}^{(0)},
\end{split}
\end{equation}
where $\phi_{j}^{(0)}$ and $\widehat{\phi}_{j}^{(0)}$ are the average values and first moments of $\phi$ at cell $I_{j}$, respectively.

\underline{\textbf{Step 2, $O(\varepsilon^{0})$ equations:}} The $O(\varepsilon^{0})$ equations can be summarized as
\begin{equation}\label{e2}
\begin{split}
&\sigma_{t,j}\left(\psi_{m,j}^{(1)}-\frac{1}{2}\sum\limits_{k=1}^{M}\psi_{k,j}^{(1)}\omega_{k}\right)
	=-\frac{\mu_{m}}{\Delta x}\left(\psi_{m,j+\frac{1}{2}}^{(0)}-\psi_{m,j-\frac{1}{2}}^{(0)}\right),\\
&\sigma_{t,j}\left(\widehat{\psi}_{m,j}^{(1)}-\frac{1}{2}\sum\limits_{k=1}^{M}\widehat{\psi}_{k,j}^{(1)}\omega_{k}\right)
	=-\frac{\mu_{m}}{2\Delta x}\left(\psi_{m,j+\frac{1}{2}}^{(0)}+\psi_{m,j-\frac{1}{2}}^{(0)}-\phi_{j}^{(0)}\right),
\end{split}
\end{equation}
where
$\psi_{m,\frac{1}{2}}^{0}=f_{m}$ if $\mu_{m}>0$, $\psi_{m,J+\frac{1}{2}}^{0}=g_{m}$ if $\mu_{m}<0$.
Following the linear HWENO reconstruction \eqref{linear21}, \eqref{linear22} and \eqref{eps0}, we have
\begin{equation}\label{eps11}
\psi_{m,j+\frac{1}{2}}^{(0)} =\frac{1}{2}\left(\frac{13}{108}\phi_{j-1}^{(0)}+\frac{7}{12}\phi_{j}^{(0)}+\frac{8}{27}\phi_{j+1}^{(0)}+\frac{25}{54}\widehat{\phi}_{j-1}^{(0)}+\frac{241}{54}\widehat{\phi_{j}}^{(0)}-\frac{28}{27}\widehat{\phi}_{j+1}^{(0)}\right),
\end{equation}
for $1\leq j \leq J$ if $\mu_{m}>0$, or
\begin{equation}\label{eps12}
\psi_{m,j+\frac{1}{2}}^{(0)}=\frac{1}{2}\left(\frac{8}{27}\phi_{j}^{(0)}+\frac{7}{12}\phi_{j+1}^{(0)}+\frac{13}{108}\phi_{j+2}^{(0)}+\frac{28}{27}\widehat{\phi}_{j}^{(0)}-\frac{241}{54}\widehat{\phi}_{j+1}^{(0)}-\frac{25}{54}\widehat{\phi}_{j+2}^{(0)}\right),
\end{equation}
$0\leq j\leq J-1$ if $\mu_{m}<0$.
Multiplying Eqs. \eqref{e2} by Gauss quadrature weights $\omega_{m}$, and summing over $m$, we find that the left
sides vanish, and the right sides yield the solvability conditions (after using \eqref{gaussweight}, more specifically, $\sum_{m}\mu_{m}\omega_{m}=0$)
\begin{equation}\label{eps00}
0=\sum_{m=1}^{M}\mu_{m}\psi_{m,j+\frac{1}{2}}^{(0)}\omega_{m},\quad 0\leq j\leq J,
\end{equation}
which must be satisfied for a solution of equation \eqref{e2} to exist.

We start by considering the case of $j=0$. The combination of \eqref{eps12} and \eqref{eps00} leads to
\begin{equation}\label{detail1}
0=\frac{1}{2}\left(\frac{8}{27}\phi_{0}^{(0)}+\frac{7}{12}\phi_{1}^{(0)}+\frac{13}{108}\phi_{2}^{(0)}+\frac{28}{27}\widehat{\phi}_{0}^{(0)}-\frac{241}{54}\widehat{\phi}_{1}^{(0)}-\frac{25}{54}\widehat{\phi}_{2}^{(0)}\right)
\sum_{\mu_{m}<0}\mu_{m}\omega_{m}+\sum_{\mu_{m}>0}\mu_{m}f_{m}\omega_{m},
\end{equation}
Let us define $\gamma$ as in \cite{larsen2}
\begin{equation*}
\gamma=2\sum_{\mu_{m}>0}\mu_{m}\omega_{m}\approx1.
\end{equation*}
Combined with \eqref{gaussweight}, we obtain
\begin{equation}\label{gamma}
\sum_{\mu_{m}>0}\mu_{m}\omega_{m}=\frac{\gamma}{2}~ \text{ and~} \sum_{\mu_{m}<0}\mu_{m}\omega_{m}=-\frac{\gamma}{2},
\end{equation}
therefore the equation \eqref{detail1} yields
\begin{equation}\label{e31}
\frac{8}{27}\phi_{0}^{(0)}+\frac{7}{12}\phi_{1}^{(0)}+\frac{13}{108}\phi_{2}^{(0)}+
\frac{28}{27}\widehat{\phi}_{0}^{(0)}-\frac{241}{54}\widehat{\phi}_{1}^{(0)}-\frac{25}{54}\widehat{\phi}_{2}^{(0)}=\frac{4}{\gamma}\sum_{\mu_{m}>0}\mu_{m}f_{m}\omega_{m}.\end{equation}
Next, for any $j$ satisfying $1\leq j\leq J-1$ or $j=J$, we follow the similar approach to combine \eqref{eps00}, \eqref{eps11}, \eqref{eps12} and \eqref{gamma} and derive
\begin{equation}\label{e32}
\begin{split}
&\frac{13}{108}\phi_{j-1}^{(0)}+\frac{7}{12}\phi_{j}^{(0)}+\frac{8}{27}\phi_{j+1}^{(0)}+\frac{25}{54}\widehat{\phi}_{j-1}^{(0)}+\frac{241}{54}\widehat{\phi}_{j}^{(0)}-\frac{28}{27}\widehat{\phi}_{j+1}^{(0)} \\
&\qquad =\frac{8}{27}\phi_{j}^{(0)}+\frac{7}{12}\phi_{j+1}^{(0)}+\frac{13}{108}\phi_{j+2}^{(0)}+\frac{28}{27}\widehat{\phi}_{j}^{(0)}-\frac{241}{54}\widehat{\phi}_{j+1}^{(0)}-\frac{25}{54}\widehat{\phi}_{j+2}^{(0)},\quad 1\leq j\leq J-1;
\end{split}
\end{equation}
\begin{equation}\label{e33}
\frac{13}{108}\phi_{J-1}^{(0)}+\frac{7}{12}\phi_{J}^{(0)}+\frac{8}{27}\phi_{J+1}^{(0)}+\frac{25}{54}\widehat{\phi}_{J-1}^{(0)}+\frac{241}{54}\widehat{\phi}_{J}^{(0)}-\frac{28}{27}\widehat{\phi}_{J+1}^{(0)}=\frac{4}{\gamma}\sum_{\mu_{m}<0}|\mu_{m}|g_{m}\omega_{m}.
\end{equation}

Let us define the following cell interface notations:
\begin{align*}
\displaystyle
&\phi_{\frac{1}{2}}^{(0)}=\frac{8}{27}\phi_{0}^{(0)}+\frac{7}{12}\phi_{1}^{(0)}+\frac{13}{108}\phi_{2}^{(0)}+
\frac{28}{27}\widehat{\phi}_{0}^{(0)}-\frac{241}{54}\widehat{\phi}_{1}^{(0)}-\frac{25}{54}\widehat{\phi}_{2}^{(0)};\\
\displaystyle
&\phi_{j+\frac{1}{2}}^{(0)}=\frac{13}{108}\phi_{j-1}^{(0)}+\frac{7}{12}\phi_{j}^{(0)}+\frac{8}{27}\phi_{j+1}^{(0)}+\frac{25}{54}\widehat{\phi}_{j-1}^{(0)}+\frac{241}{54}\widehat{\phi}_{j}^{(0)}-\frac{28}{27}\widehat{\phi}_{j+1}^{(0)}\\ \nonumber
\displaystyle
&\qquad
=\frac{8}{27}\phi_{j}^{(0)}+\frac{7}{12}\phi_{j+1}^{(0)}+\frac{13}{108}\phi_{j+2}^{(0)}+\frac{28}{27}\widehat{\phi}_{j}^{(0)}-\frac{241}{54}\widehat{\phi}_{j+1}^{(0)}-\frac{25}{54}\widehat{\phi}_{j+2}^{(0)},\quad 1\leq j\leq J-1;\\
\displaystyle
&\phi_{J+\frac{1}{2}}^{(0)}=\frac{13}{108}\phi_{J-1}^{(0)}+\frac{7}{12}\phi_{J}^{(0)}+\frac{8}{27}\phi_{J+1}^{(0)}+\frac{25}{54}\widehat{\phi}_{J-1}^{(0)}+\frac{241}{54}\widehat{\phi}_{J}^{(0)}-\frac{28}{27}\widehat{\phi}_{J+1}^{(0)},
\end{align*}
which can be denoted as
\begin{equation} \label{eq:split}
\begin{split}
\phi_{j+\frac{1}{2}}^{(0)}=&L_{1}\left(\phi_{j}^{(0)},\phi_{j+1}^{(0)},\phi_{j+2}^{(0)},\widehat{\phi}_{j}^{(0)},\widehat{\phi}_{j+1}^{(0)},\widehat{\phi}_{j+2}^{(0)}\right);\quad 0\leq j\leq J-1;\\
\phi_{j+\frac{1}{2}}^{(0)}=&L_{2}\left(\phi_{j-1}^{(0)},\phi_{j}^{(0)},\phi_{j+1}^{(0)},\widehat{\phi}_{j-1}^{(0)},\widehat{\phi}_{j}^{(0)},\widehat{\phi}_{j+1}^{(0)}\right);\quad 1\leq j\leq J;\\
\end{split}
\end{equation}
with $L_{1}$ and $L_{2}$ being two linear operators. In the ghost cells $I_0$ and $I_{J+1}$, we use the fifth order extrapolations
to evaluate their cell average and first order moments in this paper, and have
\begin{equation}\label{matrixghost}
\phi_{0}^{(0)}=\sum_{j=1}^{5}c_{j}\phi_{j}^{(0)},~~~  \phi_{J+1}^{(0)}=\sum_{j=J-4}^{J}c'_{j}\phi_{j}^{(0)},~~~
\widehat{\phi}_{0}=\sum_{j=1}^{5}\widehat{c}_{j}\widehat{\phi}_{j}^{(0)},~~~  \widehat{\phi}_{J+1}=\sum_{j=J-4}^{J}\widehat{c'}_{j}\phi_{j}^{(0)},
\end{equation}
where $c_{j},c'_{j},\widehat{c}_{j}$ and $\widehat{c'}_{j}$ are constants computed by Lagrange interpolating.
Therefore, by combining the linear relations \eqref{eq:split} and \eqref{matrixghost}, and then inverting them,
we can obtain two linear operators $L$ and $\widehat{L}$, such that
\begin{equation*}
  \phi_{j}^{(0)}=L\left(\phi_{\frac{1}{2}}^{(0)},\cdots,\phi_{J+\frac{1}{2}}^{(0)}\right), \qquad
    \widehat{\phi}_{j}^{(0)}=\widehat{L}\left(\phi_{\frac{1}{2}}^{(0)},\cdots,\phi_{J+\frac{1}{2}}^{(0)}\right),\qquad 1\leq j\leq J.
\end{equation*}

At the end of this step, let us summarize the results that are derived from Equations \eqref{e31}-\eqref{e33} and will be used later:
\begin{subequations} \label{e51-e54}
\begin{align}
\displaystyle
&\phi_{\frac{1}{2}}^{(0)}=\frac{4}{\gamma}\sum_{\mu_{m}>0}\mu_{m}f_{m}\omega_{m};\label{e51}\\
\displaystyle
&  \phi_{j}^{(0)}=L\left(\phi_{\frac{1}{2}}^{(0)},\cdots,\phi_{J+\frac{1}{2}}^{(0)}\right), \quad 1\leq j\leq J,
\label{e52}\\
\displaystyle
&\widehat{\phi}_{j}^{(0)}=\widehat{L}\left(\phi_{\frac{1}{2}}^{(0)},\cdots,\phi_{J+\frac{1}{2}}^{(0)}\right),\quad 1\leq j\leq J,\label{e53}\\
 \displaystyle
&\phi_{J+\frac{1}{2}}^{(0)}=\frac{4}{\gamma}\sum_{\mu_{m}<0}|\mu_{m}|g_{m}\omega_{m},\label{e54}
\end{align}
\end{subequations}
and
\begin{equation}
\psi_{m,j+\frac{1}{2}}^{(0)}=
\begin{cases}
f_{m},\quad &j=0,\quad \mu_{m}>0;\\
\frac{1}{2}\phi_{j+\frac{1}{2}}^{(0)},\quad &\begin{cases}
1\leq j\leq J, \quad &\mu_{m}>0;\\
0\leq j\leq J-1,\quad &\mu_{m}<0;
\end{cases}\\
g_{m},\quad &j=J,\quad \mu_{m}<0.
\end{cases}
\end{equation}
Also, note that the general solution of Equation \eqref{e2} takes the form
\begin{equation}\label{e6}
\psi_{m,j}^{(1)}=\frac{1}{2}\phi_{j}^{(1)}-\frac{\mu_{m}}{\sigma_{t,j}\Delta x}(\psi_{m,j+\frac{1}{2}}^{(0)}-\psi_{m,j-\frac{1}{2}}^{(0)}), \quad 1\leq j\leq J.\\
\end{equation}

\underline{\textbf{Step 3, $O(\varepsilon^{1})$ equations:}} Next, we consider the $O(\varepsilon^{1})$ equations, which take the form
\begin{equation}\label{e7}
\begin{split}
&\sigma_{t,j}\left(\psi_{m,j}^{(2)}-\frac{1}{2}\sum\limits_{k=1}^{M}\psi_{k,j}^{(2)}\omega_{k}\right)=-\frac{\mu_{m}}{\Delta x}\left(\psi_{m,j+\frac{1}{2}}^{(1)}-\psi_{m,j-\frac{1}{2}}^{(1)}\right)+\frac{1}{2}\left(-\sigma_{a,j}\phi_{j}^{(0)}+Q_{j}\right),\\
&\sigma_{t,j}\left(\widehat{\psi}_{m,j}^{(2)}-\frac{1}{2}\sum\limits_{k=1}^{M}\widehat{\psi}_{k,j}^{(2)}\omega_{k}\right)=-\frac{\mu_{m}}{2\Delta x}\left(\psi_{m,j+\frac{1}{2}}^{(1)}+\psi_{m,j-\frac{1}{2}}^{(1)}-2\psi_{m,j}^{(1)}\right)+\frac{1}{2}\left(-\sigma_{a,j}\widehat{\phi}_{j}^{(0)}+\widehat{Q}_{j}\right).\\
\end{split}
\end{equation}
The solvability conditions of these equations are
\begin{equation}\label{e8a}
\sum_{m=1}^{M}\mu_{m}\psi_{m,j+\frac{1}{2}}^{(1)}\omega_{m}-\sum_{m=1}^{M}\mu_{m}\psi_{m,j-\frac{1}{2}}^{(1)}\omega_{m}=\Delta x\left(-\sigma_{a,j}\phi_{j}^{(0)}+Q_{j}\right),\\
\end{equation}
\begin{equation}\label{e8b}
\sum_{m=1}^{M}\mu_{m}\psi_{m,j+\frac{1}{2}}^{(1)}\omega_{m}+\sum_{m=1}^{M}\mu_{m}\psi_{m,j-\frac{1}{2}}^{(1)}\omega_{m}=2\sum_{m=1}^{M}\mu_{m}\psi_{m,j}^{(1)}\omega_{m}
+2\Delta x\left(-\sigma_{a,j}\widehat{\phi}_{j}^{(0)}+\widehat{Q}_{j}\right).
\end{equation}
Adding Equation \eqref{e8a} over the $j$th and $(j+1)$th cells, and taking the difference of Equation \eqref{e8b} at cells $j$ and $j+1$,
yield two equalities with the same left side, which lead to the equivalence of their right sides. Therefore, we have
 \begin{align}
 \displaystyle&\sum_{m=1}^{M}\mu_{m}\left(\psi_{m,j+1}^{(1)}-\psi_{m,j}^{(1)}\right)\omega_{m}+\frac{\Delta x}{2}\left[\sigma_{a,j+1}\left(\phi_{j+1}^{(0)}-2\widehat{\phi}_{j+1}^{(0)}\right)+\sigma_{a,j}\left(\phi_{j}^{(0)}+2\widehat{\phi}_{j}^{(0)}\right)\right]	\notag \\
 \displaystyle&\hskip4cm =\frac{\Delta x}{2}\left[(Q_{j+1}-2\widehat{Q}_{j+1})+(Q_{j}+2\widehat{Q}_{j})\right],\quad 1\leq j\leq J-1.  \label{e8c}
 \end{align}

\underline{\textbf{Step 4, diffusion equation:}}
In this last step, we combine the results in the previous steps, and show that the solution $\phi_j^{(0)}$ satisfying an equation which is a consistent numerical
discretization of the diffusion equation \eqref{diffusion}.

We first plug in Equations \eqref{e51-e54}-\eqref{e6} into Equation \eqref{e8c} and obtain
 \begin{align}\label{e9a1}
\displaystyle&\hskip-1cm-\frac{1}{3\sigma_{t,j+1}\Delta x}(\phi_{j+\frac{3}{2}}^{(0)}-\phi_{j+\frac{1}{2}}^{(0)})+\frac{1}{3\sigma_{t,j}\Delta x}(\phi_{j+\frac{1}{2}}^{(0)}-\phi_{j-\frac{1}{2}}^{(0)})\nonumber \\
\displaystyle&+\frac{\Delta x}{2}\left[\sigma_{a,j+1}({\phi}_{j+1}^{(0)}-2{\widehat{{\phi}}}_{j+1}^{(0)})+\sigma_{a,j}({\phi}_{j}^{(0)}+2{\widehat{{\phi}}}_{j}^{(0)})\right]\nonumber \\
\displaystyle&=\frac{\Delta x}{2}\left[(Q_{j+1}-2\widehat{Q}_{j+1})+(Q_{j}+2\widehat{Q}_{j})\right],\quad 1\leq j\leq J-1.
 \end{align}
We consider the Taylor expansion of $\phi(x)$ and ${\phi}(x)\frac{x-x_{j}}{\Delta x}$ at any point $x_{*}\in[x_{j-\frac{1}{2}},x_{j+\frac{1}{2}}]$, which leads to
\begin{equation}\label{taylor}
  \begin{split}
     \phi(x)=&~\phi_{*}+\phi'_{*}(x-x_{*})+\frac{\phi''_{*}}{2}(x-x_{*})^2+O(\Delta x^{3}),\\
     \phi(x)\frac{x-x_{j}}{\Delta x}=&~\phi_{*}\frac{x_{*}-x_{j}}{\Delta x}+\left(\phi'_{*}\frac{x_{*}-x_{j}}{\Delta x}+\frac{\phi_{*}}{\Delta x}\right)(x-x_{*})\\
     &~+\left(\phi''_{*}\frac{x_{*}-x_{j}}{2\Delta x}+\frac{\phi'_{*}}{\Delta x}\right)(x-x_{*})^2+O(\Delta x^{3}).
   \end{split}
\end{equation}
Taking $x_{*}=x_{j-\frac{1}{2}}$, multiplying Equation \eqref{taylor} with $\frac{1}{\Delta x}$, and integrating on the cell $I_{j}$ yield
\begin{equation*}
  \begin{split}
{\phi}_{j}=&~\phi_{j-\frac{1}{2}}+\frac{\Delta x}{2}{\phi'}_{j-\frac{1}{2}}+\frac{\Delta x^{2}}{6}{\phi''}_{j-\frac{1}{2}}+O(\Delta x^3),\\
\widehat{{\phi}}_{j}=&~\frac{\Delta x}{12}{\phi'}_{j-\frac{1}{2}}-\frac{\Delta x^{2}}{12}{\phi''}_{j-\frac{1}{2}}+O(\Delta x^3),
   \end{split}
\end{equation*}
following \eqref{phi_j} and \eqref{phi_jh}, hence
\begin{equation}\label{ppart5}
{\phi}_{j}-2\widehat{{\phi}}_{j}={\phi}_{j-\frac{1}{2}}+\frac{\Delta x}{3}{\phi'}_{j-\frac{1}{2}}+\frac{\Delta x^{2}}{3}{\phi''}_{j-\frac{1}{2}}+O(\Delta x^3).
\end{equation}
Similarly, we can take $x_{*}=x_{j+\frac{1}{2}}$ and obtain
\begin{equation*}
{\phi}_{j}+2\widehat{{\phi}}_{j}={\phi}_{j+\frac{1}{2}}-\frac{\Delta x}{3}{\phi'}_{j+\frac{1}{2}}+\frac{\Delta x^{2}}{3}{\phi''}_{j+\frac{1}{2}}+O(\Delta x^3).
\end{equation*}
Hence, the equation \eqref{e9a1} becomes
\begin{align}\label{ee9aa}
\displaystyle&\hskip-1cm-\frac{1}{3\sigma_{t,j+1}\Delta x}(\phi_{j+\frac{3}{2}}^{(0)}-\phi_{j+\frac{1}{2}}^{(0)})+\frac{1}{3\sigma_{t,j}\Delta x}(\phi_{j+\frac{1}{2}}^{(0)}-\phi_{j-\frac{1}{2}}^{(0)})\nonumber \\
\displaystyle& +\frac{\Delta x}{2}(\sigma_{a,j+1}{\phi}_{j+\frac{1}{2}}^{(0)}+\sigma_{a,j}{\phi}_{j+\frac{1}{2}}^{(0)})+O(\Delta x^{2})\nonumber \\
\displaystyle&=\Delta x {Q}_{j+\frac{1}{2}}+O(\Delta x^{3}),\quad 1\leq j\leq J-1,
 \end{align}
where the right side of equation utilized the equality
\begin{equation*}
(Q_{j+1}-2\widehat{Q}_{j+1})+(Q_{j}+2\widehat{Q}_{j})=2{Q}_{j+\frac{1}{2}}+O(\Delta x^{2}),
\end{equation*}
derived in the similar way.
 Equation \eqref{ee9aa}, combined with the boundary conditions
  \begin{equation}\label{ee9bc}
 \phi_{\frac{1}{2}}^{(0)}=\frac{4}{\gamma}\sum\limits_{\mu_{m}>0}\mu_{m}f_{m}\omega_{m}, \quad \mathrm{and} \quad
  \phi_{J+\frac{1}{2}}^{(0)}=\frac{4}{\gamma}\sum\limits_{\mu_{m}<0}|\mu_{m}|g_{m}\omega_{m},
 \end{equation}
provides a consistent numerical discretization of the diffusion equation \eqref{diffusion}.

In summary, we obtain the following expression for the cell-edge angular fluxes
 \begin{equation}
 \psi_{m,j+\frac{1}{2}}=\left\{
  \begin{aligned}
 &f_{m},\quad &j=0,\mu_{m}>0\\
 &\frac{2}{\gamma}\sum_{\mu_{m}>0}\mu_{m}f_{m}\omega_{m},\quad &j=0,\mu_{m}<0\\
 &\frac{1}{2}\phi_{j+\frac{1}{2}}^{(0)},\quad &1\leq j\leq J-1\\
  &\frac{2}{\gamma}\sum_{\mu_{m}<0}|\mu_{m}|g_{m}\omega_{m},\quad &j=J, \mu_{m}>0\\
  &g_{m},\quad &j=J,\mu_{m}<0
  \end{aligned}
 \right\}+O(\varepsilon),\end{equation}
with $\phi_{\frac{1}{2}}^{(0)}$ and $\phi_{J+\frac{1}{2}}^{(0)}$ defined in \eqref{ee9bc}, and $\phi_{j+\frac{1}{2}}^{(0)}$ being the solution of \eqref{ee9aa}. For the cell-average angular fluxes, we have
\begin{equation*}
 \psi_{m,j}=\frac{1}{2}{L}\left(\phi_{\frac{1}{2}}^{(0)},\cdots,\phi_{J+\frac{1}{2}}^{(0)}\right)+O(\varepsilon),
\end{equation*}
with $L$ being the inverse operator to convert the cell interface values into the cell average values.
Therefore, when taking the limit as $\varepsilon$ approaches zero, the numerical solutions $\psi_{m,j+\frac{1}{2}}$ reduce to the solutions of the diffusion equation \eqref{diffusion} satisfying the stable and consistent method \eqref{ee9aa}. This is the AP property that is desired for the finite volume HWENO method.


\section{Multidimensional $S_{N}$ transport equations}\label{sec3}
\setcounter{equation}{0}\setcounter{figure}{0}\setcounter{table}{0}

In this section, we will discuss the finite volume HWENO FSM for multi-dimensional $S_{N}$ transport equation \eqref{Eq:transport} with isotropic scattering neutron source. Two dimensions will be used as an example to describe the HWENO FSM method, and the proposed method can be directly extended to any dimension. We consider the following two-dimensional equation
\begin{align}\label{2dsn}
&\mu\frac{\partial}{\partial x}\psi(x,y,\mu,\eta)+\eta\frac{\partial}{\partial y}\psi(x,y,\mu,\eta)+\frac{\sigma_{t,ij}}{\varepsilon}\psi(x,y,\mu,\eta) \notag \\
&\quad
=\frac{1}{4}\left(\frac{\sigma_{t,ij}}{\varepsilon}-\varepsilon\sigma_{a,ij}\right)\phi(x,y)+\frac{\varepsilon}{4}Q(x,y), \hskip10mm (x,y)\in \Omega, ~(\mu,\eta)\in[-1,1]\times[-1,1],  \\
&
\psi(x,y,\mu,\eta)=f(x,y,\mu,\eta), \hskip36mm (x,y)\in \Gamma^-,~(\mu,\eta)\in[-1,1]\times[-1,1],
\end{align}
where $\mu$ and $\eta$ represent cosine values of the angles between the neutron direction and $x$-axis and $y$-axis, respectively.
Here $\phi$ and $\psi$ are the scalar flux and angular flux with
$$\phi(x,y)=\iint_{[-1,1]\times[-1,1]}\psi(x,y,\mu,\eta)d\mu d\eta=\sum\limits_{m,n=1}^{M}\omega_{m}\omega_{n}\psi(x,y,\mu_{m},\eta_{n}),$$
where $\omega_{m}$ and $\omega_{n}$ are the level symmetric quadrature weights.

Assume the computational domain has been divided into cells $I_{i,j}=J_i\times K_j=[x_{i-\frac{1}{2}},x_{i+\frac{1}{2}}]\times[y_{j-\frac{1}{2}},y_{j+\frac{1}{2}}]$, 
with $i=1,\cdots,N_{x}$, $j=1,\cdots,N_{y}$. We denote the cell center as $(x_{i},y_{j})=((x_{i-\frac{1}{2}}+x_{i+\frac{1}{2}})/2,(y_{j-\frac{1}{2}}+y_{j+\frac{1}{2}})/2)$ and the mesh siz as $\Delta x_{i}=x_{i+\frac{1}{2}}-x_{i-\frac{1}{2}}$, $\Delta y_{j}=y_{j+\frac{1}{2}}-y_{j-\frac{1}{2}}$. The cell average and first order moments of the unknown are denoted as
\begin{equation}
\begin{split}
&\psi_{i,j}^{m,n}=\frac{1}{\Delta x_{i} \Delta y_{j}}\iint_{I_{i,j}}\psi(x,y,\mu_{m},\eta_{n})dxdy, \\
&\widehat{\psi}_{i,j}^{m,n}=\frac{1}{\Delta x_{i} \Delta y_{j}}\iint_{I_{i,j}}\psi(x,y,\mu_{m},\eta_{n})\frac{x-x_{i}}{\Delta x_{i}}dxdy,\\
&\widetilde{\psi} _{i,j}^{m,n}=\frac{1}{\Delta x_{i} \Delta y_{j}}\iint_{I_{i,j}}\psi(x,y,\mu_{m},\eta_{n})\frac{y-y_{j}}{\Delta y_{j}}dxdy, \\
&\widehat{\widetilde{\psi} }_{i,j}^{m,n}=\frac{1}{\Delta x_{i} \Delta y_{j}}\iint_{I_{i,j}}\psi(x,y,\mu_{m},\eta_{n})\frac{x-x_{i}}{\Delta x_{i}}\frac{y-y_{j}}{\Delta y_{j}}dxdy.
\end{split}
\end{equation}
For simplicity, we ignore the superscript $m,n$ without causing any confusion. Similarly, we can define the cell average and moments of $\phi$ as
$\phi_{i,j},\widehat{\phi}_{i,j},\widetilde{\phi} _{i,j},\widehat{\widetilde{\phi} }_{i,j}$, and those of $Q$ as
$Q_{i,j},\widehat{Q}_{i,j},\widetilde{Q} _{i,j},\widehat{\widetilde{Q} }_{i,j}$.
 We multiply \eqref{2dsn} by $\frac{1}{\Delta x_{i} \Delta y_{j}}$, $\frac{x-x_{i}}{\Delta x_{i}^2 \Delta y_{j}}$, $\frac{y-y_{j}}{\Delta x_{i}\Delta y_{j}^2}$ and $\frac{x-x_{i}}{\Delta x_{i}^2}\frac{y-y_{j}}{\Delta y_{j}^2}$, respectively, and then integrate them on the cell $I_{i,j}$. Applying integration by parts and replacing the cell interface values by the Godunov numerical flux as discussed in one-dimensional setting in Section \ref{sec:hweno1d}, we have the following HWENO discretization
\begin{align}
&\frac{\mu_m}{\Delta x_{i}\Delta y_{j}}\int_{K_{j}}[\psi^\pm(x_{i+\frac{1}{2}},y,\mu_{m},\eta_{n})-\psi^\pm(x_{i-\frac{1}{2}},y,\mu_{m},\eta_{n})]dy  \notag \\
&\hskip25mm
+\frac{\eta_n}{\Delta x_{i}\Delta y_{j}}\int_{J_{i}}[\psi^\pm(x,y_{j+\frac{1}{2}},\mu_{m},\eta_{n})-\psi^\pm(x,y_{j-\frac{1}{2}},\mu_{m},\eta_{n})]dx \notag \\
&\hskip25mm
+\frac{\sigma_{t,ij}}{\varepsilon}\psi_{i,j}^{m,n}
=\frac{1}{4}\left(\frac{\sigma_{t,ij}}{\varepsilon}-\varepsilon\sigma_{a,ij}\right)\sum\limits_{k,l=1}^{M}\psi_{i,j}^{k,l}\omega_{k}\omega_{l}+\frac{\varepsilon}{4}Q_{i,j},\label{discre2d1}\\
&\frac{\mu_m}{2\Delta x_{i}\Delta y_{j}}\int_{K_{j}}[\psi^\pm(x_{i+\frac{1}{2}},y,\mu_{m},\eta_{n})+\psi^\pm(x_{i-\frac{1}{2}},y,\mu_{m},\eta_{n})]dy-\frac{\mu_m}{\Delta x_{i}}\psi_{i,j}^{m,n} \notag \\
&\hskip25mm
+\frac{\eta_n}{\Delta x_{i}\Delta y_{j}}\int_{J_{i}}\frac{x-x_i}{\Delta x_{i}}[\psi^\pm(x,y_{j+\frac{1}{2}},\mu_{m},\eta_{n})-\psi^\pm(x,y_{j-\frac{1}{2}},\mu_{m},\eta_{n})]dx\notag \\
&\hskip25mm
+\frac{\sigma_{t,ij}}{\varepsilon}\widehat{\psi}_{i,j}^{m,n}=\frac{1}{4}\left(\frac{\sigma_{t,ij}}{\varepsilon}-\varepsilon\sigma_{a,ij}\right)\sum\limits_{k,l=1}^{M}\widehat{\psi}_{i,j}^{k,l}\omega_{k}\omega_{l}+\frac{\varepsilon}{4}\widehat{Q}_{i,j},\label{discre2d2}\\
&\frac{\mu_m}{\Delta x_{i}\Delta y_{j}}\int_{K_{j}}\frac{y-y_{j}}{\Delta y_{j}}[\psi^\pm(x_{i+\frac{1}{2}},y,\mu_{m},\eta_{n})-\psi^\pm(x_{i-\frac{1}{2}},y,\mu_{m},\eta_{n})]dy\notag \\
&\hskip25mm
+\frac{\eta_n}{2\Delta x_{i}\Delta y_{j}}\int_{J_{i}}[\psi^\pm(x,y_{j+\frac{1}{2}},\mu_{m},\eta_{n})+\psi^\pm(x,y_{j-\frac{1}{2}},\mu_{m},\eta_{n})]dx-\frac{\eta_n}{\Delta y_{j}}\psi_{i,j}^{m,n}\notag \\
&\hskip25mm
+\frac{\sigma_{t,ij}}{\varepsilon}\widetilde{\psi} _{i,j}^{m,n}=\frac{1}{4}\left(\frac{\sigma_{t,ij}}{\varepsilon}-\varepsilon\sigma_{a,ij}\right)\sum\limits_{k,l=1}^{M}\widetilde{\psi} _{i,j}^{k,l}\omega_{k}\omega_{l}+\frac{\varepsilon}{4}\widetilde{Q} _{i,j}, \label{discre2d3}\\
	\displaystyle
&\frac{\mu_m}{2\Delta x_{i}\Delta y_{j}}\int_{K_{j}}\frac{y-y_{j}}{\Delta y_{j}}[\psi^\pm(x_{i+\frac{1}{2}},y,\mu_{m},\eta_{n})+\psi^\pm(x_{i-\frac{1}{2}},y,\mu_{m},\eta_{n})]dy-\frac{\mu_m}{\Delta x_{i}}\widetilde{\psi}_{i,j}^{m,n} \notag \\
&\hskip25mm
+\frac{\eta_n}{2\Delta x_{i}\Delta y_{j}}\int_{J_{i}}\frac{x-x_{i}}{\Delta x_{i}}[\psi^\pm(x,y_{j+\frac{1}{2}},\mu_{m},\eta_{n})+\psi^\pm(x,y_{j-\frac{1}{2}},\mu_{m},\eta_{n})]dx -\frac{\eta_n}{\Delta y_{j}}\widehat{\psi}_{i,j}^{m,n} \notag \\
&\hskip25mm
+\frac{\sigma_{t,ij}}{\varepsilon}\widehat{\widetilde{\psi} }_{i,j}^{m,n}=\frac{1}{4}\left(\frac{\sigma_{t,ij}}{\varepsilon}-\varepsilon\sigma_{a,ij}\right)\sum\limits_{k,l=1}^{M}\widehat{\widetilde{\psi} }_{i,j}^{k,l}\omega_{k}\omega_{l}+\frac{\varepsilon}{4}\widehat{\widetilde{Q} }_{i,j}, \label{discre2d4}
\end{align}
where the numerical flux $\psi^\pm(x_{i+\frac{1}{2}},y,\mu_m,\eta_n)$ is chosen to $\psi^{+}$ when $\mu_m>0$ and $\psi^-$ otherwise. Similarly, the numerical flux $\psi^\pm(x,y_{j+\frac{1}{2}},\mu_m,\eta_n)$ is chosen to $\psi^{+}$ when $\eta_n>0$ and $\psi^-$ otherwise. The integrals of the flux over $K_j$ or $J_i$ are evaluated via the HWENO reconstruction to be discussed in the following subsection.

\subsection{HWENO reconstruction in 2D} \label{2dhweno}
We can use the dimension-by-dimension strategy to reconstruct these integrals in the HWENO method \eqref{discre2d1}-\eqref{discre2d4}. The procedure of these reconstructions is sketched as follows. Again, for ease of presentation, we assume the uniform mesh with $\Delta x_{i}=\Delta x$, $\Delta y_{j}=\Delta y$ in the description. We denote $\psi^{m,n}(x,y)=\psi(x,y,\mu_m,\eta_n)$ and ignore the superscript $(m,n)$ below without causing any confusion.
\begin{itemize}
  \item In the $x$-direction, we perform the one-dimensional HWENO reconstruction which was described in subsection \ref{sec:hweno1d}. Therefore, from $\{\psi_{l,j},\widehat{\psi}_{l,j}\}_{l=i-1}^{i+1}$, we can obtain $\frac{1}{\Delta y}\int_{K_{j}}\psi^\pm(x_{i\pm\frac{1}{2}},y)dy$, which is the point value in the $x$-direction and the cell-average in the $y$-direction. Similarly, we can use the values $\{\widetilde{\psi} _{l,j},\widehat{\widetilde{\psi} }_{l,j}\}_{l=i-1}^{i+1}$ to reconstruct $\frac{1}{\Delta y}\int_{K_{j}}\psi^\pm(x_{i\pm\frac{1}{2}},y)\frac{y-y_{j}}{\Delta y}dy$. Note that either ``+'' or ``-'' sign is taken, depending on whether $\mu_m$ is negative or positive.

  \item In the $y$-direction, we perform the one-dimensional HWENO reconstruction which was described in subsection \ref{sec:hweno1d}. Therefore, from $\{\psi_{i,l},\widetilde{\psi} _{i,l}\}_{l=j-1}^{j+1}$, we can obtain $\frac{1}{\Delta x}\int_{J_{i}}\psi(x,y_{j\pm\frac{1}{2}})dx$, which is the point value in the $y$-direction and the cell-average in the $x$-direction. Similarly, we can use the values $\{\widehat{\psi}_{i,l},\widehat{\widetilde{\psi} }_{i,l}\}_{l=j-1}^{j+1}$ to reconstruct $\frac{1}{\Delta x}\int_{J_{i}}\psi(x,y_{j\pm\frac{1}{2}})\frac{x-x_{i}}{\Delta x}dx$. Note that either ``+'' or ``-'' sign is taken, depending on whether $\eta_n$ is negative or positive.

\end{itemize}

\subsection{Fast sweeping idea to solve global linear system in 2D}
The proposed two-dimensional HWENO scheme for the linear transport equation \eqref{2dsn} takes the form of \eqref{discre2d1}-\eqref{discre2d4}, combined with the HWENO reconstruction to evaluate the fluxes. This is a large system involving the flux term (coupling in $x$ and $y$ directions) on the left side and the summation term (coupling in $\mu$ and $\eta$ directions) on the right side. As in one-dimensional case, the fast sweeping idea is adopted to solve this system efficiently. Let us first denote the right-hand side term of the four equations in \eqref{discre2d1}-\eqref{discre2d4} as $S_{i,j}$, $\widehat{S}_{i,j}$, $\widetilde{S}_{i,j}$ and $\widehat{\widetilde{S}}_{i,j}$, respectively.
We summarize the flowchart of HWENO FSM for the $S_{N}$ equations in 2D as follows.

\noindent \textbf{Step 1}. \emph{Initialization}:
We take $0$ as the initial value of the unknowns: $\psi_{i,j}^{m,n}$, $\widehat{\psi}_{i,j}^{m,n}$ $\widetilde{\psi}_{i,j}^{m,n}$ and $\widehat{\widetilde{\psi}}_{i,j}^{m,n}$ for all $m$, $n$, $i$ and $j$. Then we can evaluate $S_{i,j}$, $\widehat{S}_{i,j}$, $\widetilde{S}_{i,j}$ and $\widehat{\widetilde{S}}_{i,j}$.

\noindent\textbf{Step 2}. \emph{Gauss-Seidel iteration with alternating sweep}.
 We sweep the whole domain with the following four alternating orderings repeatedly for all $m$ and $n$:

(I) $i=1\rightarrow N_{x}$, $j=1\rightarrow N_{y}$: if $\mu_{m}>0$ and $\eta_{n}>0$, solve the system \eqref{discre2d1}-\eqref{discre2d4} with the appropriate boundary conditions with this order of $i,~j$. After updating the approximations $\psi_{i,j}^{m,n}$, $\widehat{\psi}_{i,j}^{m,n}$, $\widetilde{\psi}_{i,j}^{m,n}$ and $\widehat{\widetilde{\psi}}_{i,j}^{m,n}$ in the cell $I_{i,j}$, we can apply HWENO reconstruction to obtain the cell-edge flux fluxes $\frac{1}{\Delta y}\int_{K_{j}}\psi^-(x_{i+\frac{1}{2}},y)dy$ and $\frac{1}{\Delta x}\int_{J_{i}}\psi^-(x,y_{j+\frac{1}{2}})dx$ based on the most updated values of the unknowns.

(II) $i=N_{x}\rightarrow 1$, $j=1 \rightarrow N_{y}$: if $\mu_{m}<0$ and $\eta_{n}>0$, solve the system \eqref{discre2d1}-\eqref{discre2d4} with the appropriate boundary conditions with this order of $i,~j$. After updating the approximations $\psi_{i,j}^{m,n}$, $\widehat{\psi}_{i,j}^{m,n}$, $\widetilde{\psi}_{i,j}^{m,n}$ and $\widehat{\widetilde{\psi}}_{i,j}^{m,n}$ in the cell $I_{i,j}$, we can apply HWENO reconstruction to obtain the cell-edge flux fluxes $\frac{1}{\Delta y}\int_{K_{j}}\psi^+(x_{i-\frac{1}{2}},y)dy$ and $\frac{1}{\Delta x}\int_{J_{i}}\psi^-(x,y_{j+\frac{1}{2}})dx$ based on the most updated values of the unknowns.

(III) $i=1\rightarrow N_{x}$, $j=N_{y}\rightarrow 1$: if $\mu_{m}>0$ and $\eta_{n}<0$, solve the system \eqref{discre2d1}-\eqref{discre2d4} with the appropriate boundary conditions with this order of $i,~j$. After updating the approximations $\psi_{i,j}^{m,n}$, $\widehat{\psi}_{i,j}^{m,n}$, $\widetilde{\psi}_{i,j}^{m,n}$ and $\widehat{\widetilde{\psi}}_{i,j}^{m,n}$ in the cell $I_{i,j}$, we can apply HWENO reconstruction to obtain the cell-edge flux fluxes $\frac{1}{\Delta y}\int_{K_{j}}\psi^-(x_{i+\frac{1}{2}},y)dy$ and $\frac{1}{\Delta x}\int_{J_{i}}\psi^+(x,y_{j-\frac{1}{2}})dx$ based on the most updated values of the unknowns.

(IV) $i=N_{x} \rightarrow 1$, $j=N_{y}\rightarrow 1$: if $\mu_{m}<0$ and $\eta_{n}<0$, solve the system \eqref{discre2d1}-\eqref{discre2d4} with the appropriate boundary conditions with this order of $i,~j$. After updating the approximations $\psi_{i,j}^{m,n}$, $\widehat{\psi}_{i,j}^{m,n}$, $\widetilde{\psi}_{i,j}^{m,n}$ and $\widehat{\widetilde{\psi}}_{i,j}^{m,n}$ in the cell $I_{i,j}$, we can apply HWENO reconstruction to obtain the cell-edge flux fluxes $\frac{1}{\Delta y}\int_{K_{j}}\psi^+(x_{i-\frac{1}{2}},y)dy$ and $\frac{1}{\Delta x}\int_{J_{i}}\psi^+(x,y_{j-\frac{1}{2}})dx$ based on the most updated values of the unknowns.

At the boundary of the computational domain, high order extrapolations are used to compute the values at the ghost cells, which are needed for the HWENO reconstruction near the boundary. After repeating this process for all $m, n$ directions, we can compute scalar flux $\phi$ from $\psi^{m,n}$ via Gauss quadrature, and update $S_{i,j}$ etc. This completes one Gauss-Seidel iteration.

\noindent \textbf{Step 3}. \emph{Convergence}: Repeat the Gauss-Seidel iteration until the convergence criteria is satisfied. In this paper, if the scalar flux satisfies
$$\delta=||\phi^{new}-\phi^{old}||_{L_{1}}<10^{-14},$$
for two consecutive iteration steps, we stop the iteration.

\begin{remark}\label{remark2}
 The numerical tests show that HWENO FSM will not converge to machine epsilon, i.e. $\delta$ will not decrease to $10^{-14}$ for 2D examples. This has also been observed in the application of HWENO FSM method to other system \cite{hwenofsm}. To fix this, we propose to update the solution by
$$\psi_{m,j}^{new}=\omega\psi_{i,j}^{new}+(1-\omega)\psi_{i,j}^{old}, ~0<\omega<1,$$
which is shown to yield good convergence, although it may slightly increase the number of iterations. Numerically, one observes that $\omega=0.85$ is the optimal choice and will be used in our 2D numerical examples.
\end{remark}

\section{Numerical results}\label{sec-num}
\setcounter{equation}{0}\setcounter{figure}{0}\setcounter{table}{0}

We present extensive one-dimensional and two-dimensional numerical results on different model
problems, to demonstrate the diffusion limit and order of accuracy of the proposed HWENO fast sweeping method in the finite volume framework.
In all the numerical examples, $\widetilde{\varepsilon}$ in \eqref{nonweight} is taken as $10^{-6}$ unless otherwise specified. We use ``iter'' to indicate the number of iterations (noting that one iteration includes two alternating sweeping for 1D problem and four alternating sweepings for 2D problem) in all the tables. The number of grid points is assumed to be $N_x=N_y=N$ for 2D examples. All the computations are implemented by using MATLAB 2020a on ThinkPad computer with 1.80 GHz Intel Core i7 processor and 16GB RAM.

\subsection{One-dimensional problem with vacuum boundary}

\begin{figure}
\begin{center}
  \includegraphics[width=8cm]{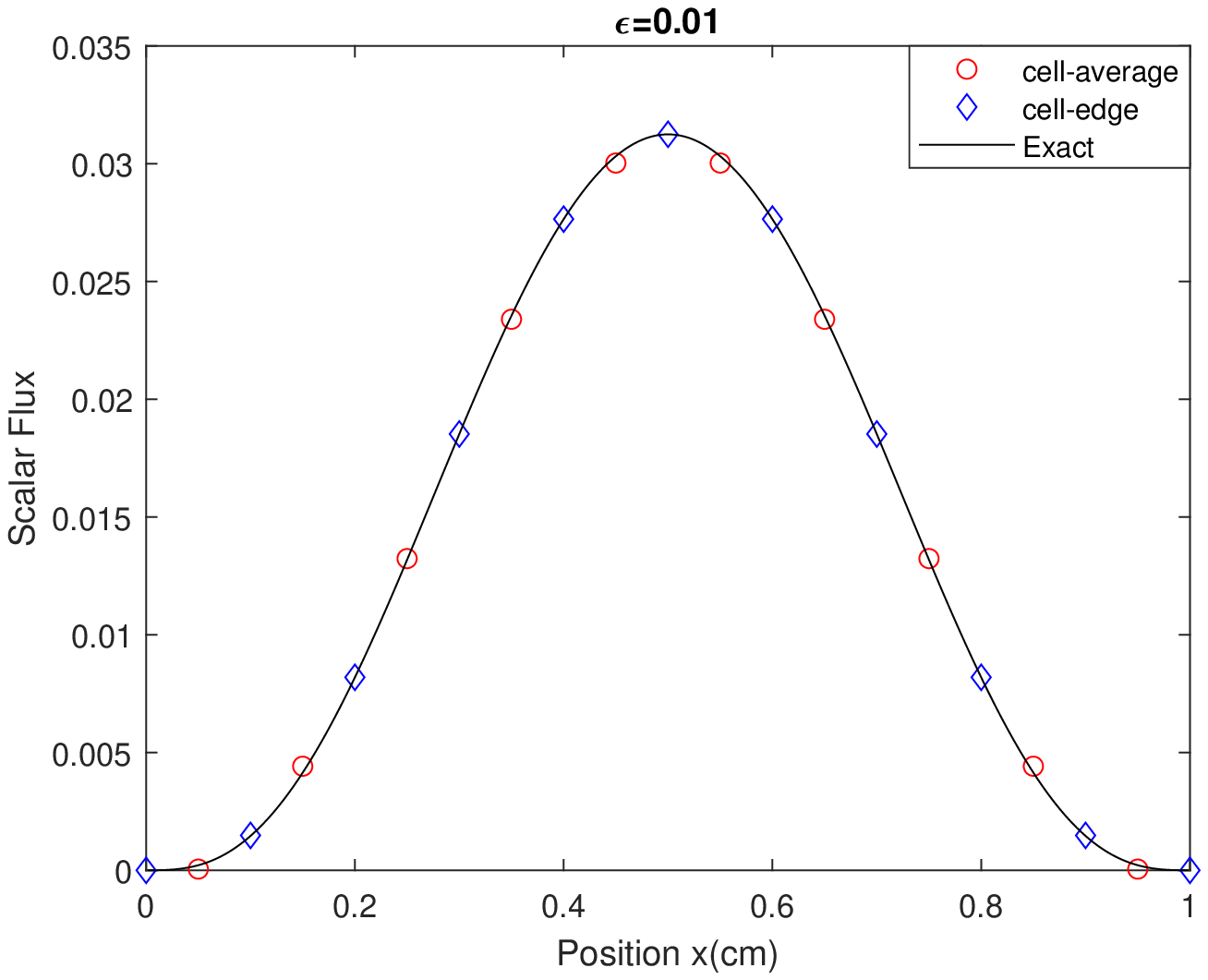}
  \includegraphics[width=8cm]{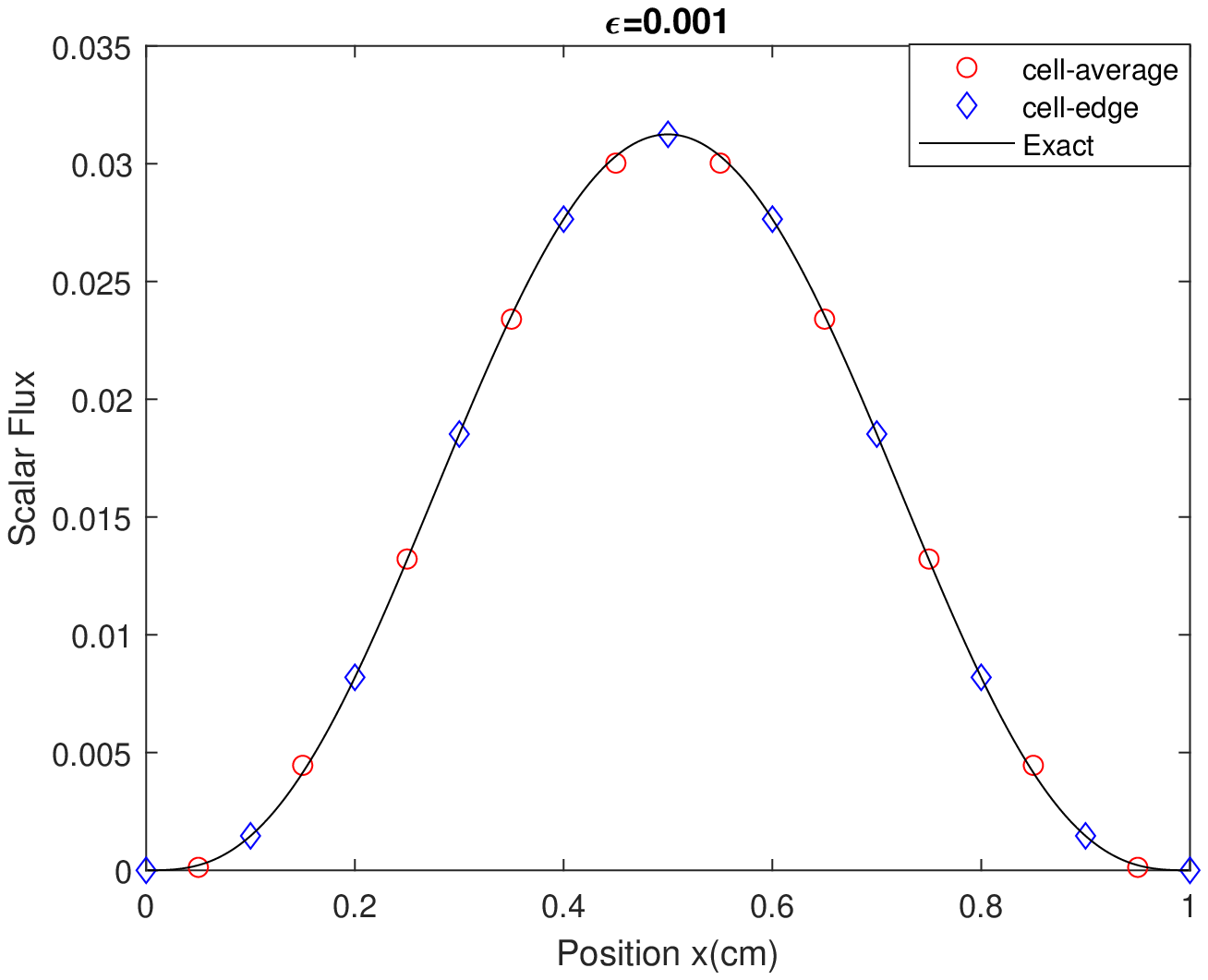}
\end{center}
 \caption{Numerical solution of Example 1 with $\Delta x=0.1$. $\varepsilon=0.01$ (left) and $\varepsilon=0.001$ (right)}\label{e2fig}
\end{figure}

\noindent \textbf{Example 1} (Accuracy test with manufactured solution).
In the first example, we consider a slab with the vacuum
boundary on both sides to test the accuracy of the proposed HWENO method.
The specifications of the problem are given as
$$L=1,~\sigma_{t}=1,~\sigma_{a}=0.8, ~ Q=\frac{2}{\varepsilon}[(3x^2-12x^3+15x^4-6x^5)\mu_{m}]+2\sigma_{a}x^{3}(1-x)^3,$$
where $L$ is the slab thickness. The manufactured exact solution of the linear transport equation is given by \cite{wangictt}
$$\psi(x,\mu_{m})=x^{3}(1-x)^3.$$
The Gauss-Legendre $S_{12}$ quadrature set is used in the angular discretization.
We have run the simulations for various choices of $\eps$.
In Table \ref{e2hytable}, we show the numerical errors, the corresponding order of accuracy and computational time of HWENO method with the hybrid strategy discussed in Remark \ref{remark2.1}. Here we only report the errors of cell-average to save space, and similar behavior has been observed for the first order moment. We can
observe that the expected high order accuracy has been observed for all choices of $\eps$. The numerical solutions with spatial size $\Delta x=0.1$, $\eps=0.01$ and $0.001$, compared with the corresponding reference solutions, are plotted in Fig. \ref{e2fig}, from which we can observe a good match of the numerical solution even for small $\eps$.

In addition, we also present the numerical errors and the corresponding order of accuracy of HWENO without the hybrid strategy in Table \ref{e2nohytable}.
One can see that the error and order are similar, while the main difference lies in the computational time. We have summarized the comparison of computational
time of the HWENO method with or without the hybrid strategy in Table \ref{timetable}. Note that the recorded time listed in Table \ref{timetable} is the total computational time of all the simulations with $N=10$, $20$, $\cdots$, $160$ and $\varepsilon=0.01$. For this example, it shows that hybrid strategy can save $70\%$ of CPU time.

\begin{table}
\caption{Example 1. The errors, order of accuracy and CPU time of HWENO method with the hybrid strategy}\label{e2hytable}
	\begin{center}
		\begin{tabular}{|c|c|c|c|c|c|c|}
			\hline
\multicolumn{7}{|c|}{$\varepsilon=1$}\\ \hline
$N$ & $L_{1}$~error&order&  $L_{\infty}$~error&order&iter &time (sec)\\\hline
10 &1.26e-05 &- &5.53e-05 &- &60 &0.006\\ \hline
20 &1.84e-07 &6.09 &1.525376e-06 &5.18 &56 &0.002\\ \hline
40 &1.96e-09 &6.55 &3.131162e-08 &5.60 &52 &0.003\\ \hline
80 &1.80e-11 &6.76 &5.617674e-10 &5.800 &49 &0.009\\ \hline
160 &1.69e-13 &6.73 &9.415637e-12 &5.89 &47 &0.013\\ \hline
\multicolumn{7}{|c|}{$\varepsilon=0.1$}\\ \hline
$N$ & $L_{1}$~error&order&  $L_{\infty}$~error&order&iter &time (sec)\\\hline
10 &5.59e-05 &- &2.59e-04 &- &886 &0.02\\ \hline
20 &9.91e-07 &5.81 &8.23e-06 &4.97 &884 &0.04\\ \hline
40 &1.40e-08 &6.14 &1.98e-07 &5.37 &883 &0.04\\ \hline
80 &1.80e-10 &6.28 &4.17e-09 &5.57 &882 &0.05\\ \hline
160 &2.06e-12 &6.44 &7.91e-11 &5.72 &882 &0.13\\ \hline
\multicolumn{7}{|c|}{$\varepsilon=0.01$}\\ \hline
$N$ & $L_{1}$~error&order&  $L_{\infty}$~error&order&iter &time (sec)\\\hline
10 &7.94e-05 &- &3.40e-04 &- &60976 &1.41\\ \hline
20 &1.68e-06 &5.55 &1.57e-05 &4.43 &60976 &1.76\\ \hline
40 &2.98e-08 &5.82 &5.36e-07 &4.87 &60969 &2.81\\ \hline
80 &5.88e-10 &5.66 &1.48e-08 &5.17 &60973 &4.31\\ \hline
160 &1.25e-11 &5.54 &3.55e-10 &5.38 &60970 &7.15\\ \hline
		\end{tabular}
	\end{center}
\end{table}
\begin{table}
\caption{Example 1. The errors, order of accuracy and CPU time of HWENO method without the hybrid strategy}\label{e2nohytable}
	\begin{center}
		\begin{tabular}{|c|c|c|c|c|c|c|}
			\hline
\multicolumn{7}{|c|}{$\varepsilon=1$}\\ \hline
$N$ & $L_{1}$~error&order&  $L_{\infty}$~error&order&iter &time (sec)\\\hline
10 &1.26e-05 & -            &5.53e-05 &-             &60 &0.006\\ \hline
20 &1.84e-07 &6.09 &1.52e-06 &5.18 &56 &0.004\\ \hline
40 &1.96e-09 &6.55 &3.13e-08 &5.60 &52 &0.006\\ \hline
80 &1.80e-11 &6.76 &5.61e-10 &5.80 &49 &0.012\\ \hline
160 &1.69e-13 &6.73 &9.41e-12 &5.89 &47 &0.022\\ \hline
\multicolumn{7}{|c|}{$\varepsilon=0.1$}\\ \hline
$N$ & $L_{1}$~error&order&  $L_{\infty}$~error&order&iter &time (sec)\\\hline
10 &5.59e-05 &-           &2.59e-04 &-             &886 &0.04\\ \hline
20 &9.91e-07 &5.81 &8.23e-06 &4.97 &884 &0.06\\ \hline
40 &1.40e-08 &6.14 &1.98e-07 &5.37 &883 &0.10\\ \hline
80 &1.80e-10 &6.28 &4.17e-09 &5.57 &882 &0.23\\ \hline
160 &2.06e-12 &6.44 &7.91e-11 &5.72 &882 &0.39\\ \hline
\multicolumn{7}{|c|}{$\varepsilon=0.01$}\\ \hline
$N$ & $L_{1}$~error&order&  $L_{\infty}$~error&order&iter &time (sec)\\\hline
10 &7.94e-05 &-            &3.40e-04 &-            &60990 &2.94\\ \hline
20 &1.68e-06 &5.55 &1.57e-05 &4.43 &60961 &4.25\\ \hline
40 &2.98e-08 &5.82 &5.36e-07 &4.87 &60969 &7.89\\ \hline
80 &5.88e-10 &5.66 &1.48e-08 &5.17 &60970 &16.38\\ \hline
160 &1.25e-11 &5.54 &3.55e-10 &5.38 &60973 &31.24\\ \hline
		\end{tabular}
	\end{center}
\end{table}

\begin{figure}
\begin{center}
  \includegraphics[width=8cm]{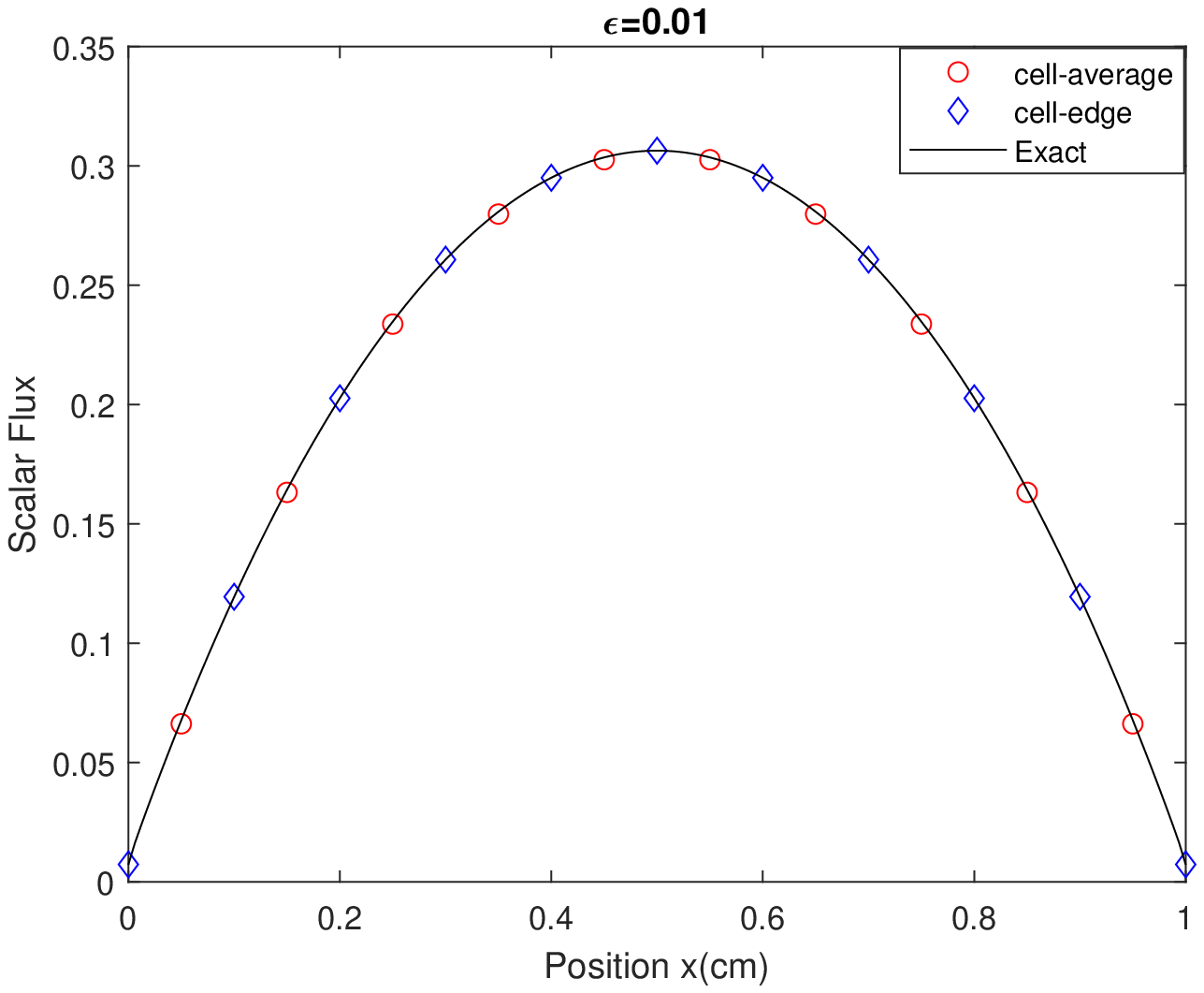}
  \includegraphics[width=8cm]{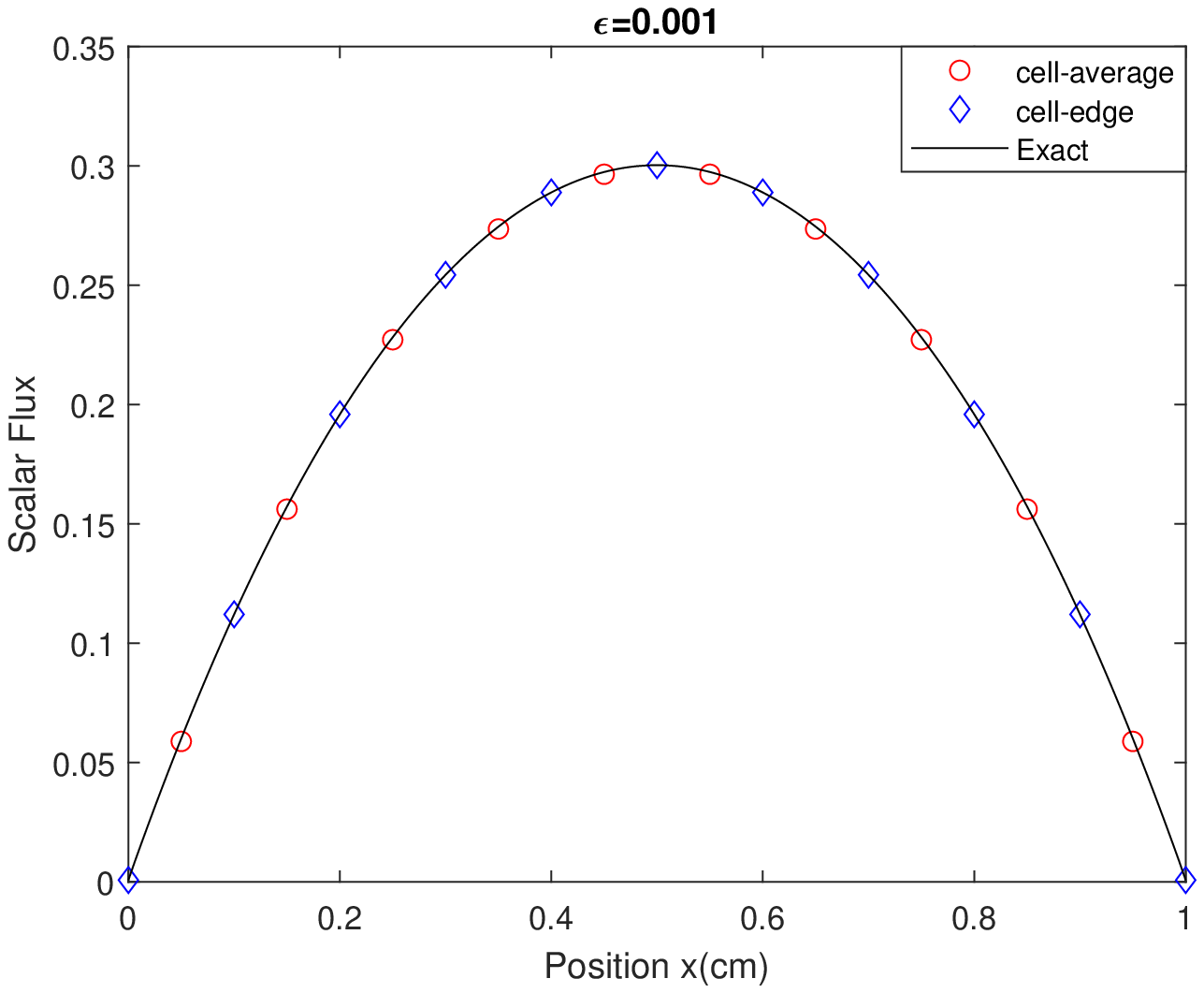}
\end{center}
 \caption{Numerical solution of Example 2 with $\Delta x=0.1$. $\varepsilon=0.01$ (left) and $\varepsilon=0.001$ (right)}\label{e1fig}
\end{figure}

\noindent \textbf{Example 2}. A slab with the vacuum
boundary on both sides with
$$L=1,~~~\sigma_{t}=1,~~~\sigma_{a}=0.8, ~~~ Q=1,$$
is considered. Again, the Gauss-Legendre $S_{12}$ quadrature set is used in the angular discretization.
The analytical solution can be obtained following the approach discussed in \cite{wanganalytic}.
We have run the simulations for various choices of $\eps$.
In Table \ref{e1hytable}, we show the numerical errors, the corresponding order of accuracy and computational time of HWENO method with the hybrid strategy discussed in Remark \ref{remark2.1}. Here we only report the errors of cell-average to save space, and similar behavior has been observed for the first order moment. We can
observe that the expected high order accuracy has been observed for all choices of $\eps$. The numerical solutions with spatial size $\Delta x=0.1$, $\eps=0.01$ and $0.001$, compared with the corresponding reference solutions, are plotted in Fig. \ref{e1fig}, from which we can observe a good match of the numerical solution even for small $\eps$.
From the computational time comparison of HWENO method with or without the hybrid strategy in Table \ref{timetable}, we observe a $80\%$ saving of CPU time when the hybrid strategy is used.

As $\varepsilon$ decreases the problem becomes thick and diffusive, and its asymptotic solution should be the same as the solution of the corresponding
diffusion equation. In Table \ref{epsilontable}, we list the errors between numerical solution of $S_{N}$ equation with different $\varepsilon$ and the exact solution of the diffusion equation, from which we can observe the error decays at the expected rate of $O(\eps)$.

\begin{table}
\caption{Example 2. The errors, order of accuracy and CPU time of HWENO method with the hybrid strategy}\label{e1hytable}
	\begin{center}
		\begin{tabular}{|c|c|c|c|c|c|c|}
			\hline
\multicolumn{7}{|c|}{$\varepsilon=1$}\\ \hline
$N$ & $L_{1}$~error&order&  $L_{\infty}$~error&order&iter &time (sec)\\\hline
10 &1.80e-05 &-            &4.60e-05 &-            &53 &0.002000\\ \hline
20 &4.07e-07 &5.47 &2.19e-06 &4.38 &50 &0.002045\\ \hline
40 &5.50e-09 &6.20 &6.03e-08 &5.18 &46 &0.002242\\ \hline
80 &5.69e-11 &6.59 &1.25e-09 &5.58 &43 &0.002936\\ \hline
160 &5.13e-13 &6.79 &2.27e-11 &5.78 &41 &0.007173\\ \hline
\multicolumn{7}{|c|}{$\varepsilon=0.1$}\\ \hline
$N$ & $L_{1}$~error&order&  $L_{\infty}$~error&order&iter &time (sec)\\ \hline
10 &2.18e-04 &-              &8.09e-04 &-           &882 &0.023520\\ \hline
20 &6.13e-05 &1.82 &3.46e-04 &1.22 &882 &0.024328\\ \hline
40 &8.90e-06 &2.78 &8.64e-05 &2.00 &883 &0.030132\\ \hline
80 &5.24e-07 &4.08 &9.11e-06 &3.24 &883 &0.047054\\ \hline
160 &1.39e-08 &5.22 &4.41e-07 &4.36 &883 &0.080591\\ \hline
\multicolumn{7}{|c|}{$\varepsilon=0.01$}\\ \hline
$N$ & $L_{1}$~error&order&  $L_{\infty}$~error&order&iter &time (sec)\\ \hline
10 &4.62e-05 &-             &1.71e-04 &-             &61034 &1.6279\\ \hline
20 &3.17e-05 &0.54 &1.78e-04 &-0.05 &61051 &1.5876\\ \hline
40 &1.89e-05 &0.74 &1.48e-04 &0.26 &61057 &2.2305\\ \hline
80 &9.46e-06 &1.00 &9.08e-05 &0.70 &61059 &3.2170\\ \hline
160 &3.57e-06 &1.40 &4.29e-05 &1.08 &61062 &5.6801\\ \hline
320 &7.85e-07 &2.18 &1.33e-05 &1.68 &61063 &10.6818\\ \hline
640 &7.26e-08 &3.430 &1.87e-06 &2.83 &61066 &20.0900\\ \hline
1280 &2.76e-09 &4.71 &1.13e-07 &4.04 &61068 &38.8375\\ \hline
		\end{tabular}
	\end{center}
\end{table}
%

\begin{table}
\caption{Example 2. The errors and order of accuracy between numerical solution of $S_{N}$ equation with different $\varepsilon$ and the exact solution of the limit diffusion equation }\label{epsilontable}
	\begin{center}
		\begin{tabular}{|c|c|c|c|c|}
			\hline
\multicolumn{5}{|c|}{$N=10$}\\ \hline
$\varepsilon$ & $L_{1}$~error&order&  $L_{\infty}$~error&order\\\hline
1 &4.89e-01 &-              &5.05e-01 &- \\ \hline
0.1 &7.16e-02 &0.83 &7.98e-02 &0.80 \\ \hline
0.01 &7.32e-03 &0.99 &8.19e-03 &0.98 \\ \hline
0.001 &7.31e-04 &1.00 &8.11e-04 &1.00 \\ \hline
0.0001 &7.31e-05 &1.00 &8.09e-05 &1.00 \\ \hline
\multicolumn{5}{|c|}{$N=20$}\\ \hline
$\varepsilon$& $L_{1}$~error&order&  $L_{\infty}$~error&order\\\hline
1 &4.89e-01 &-       &5.05e-01 &- \\ \hline
0.1 &7.20e-02 &0.83 &8.05e-02 &0.79 \\ \hline
0.01 &7.38e-03 &0.98 &8.60e-03 &0.97 \\ \hline
0.001 &7.37e-04 &1.00 &8.49e-04 &1.00 \\ \hline
0.0001 &7.37e-05 &1.00 &8.46e-05 &1.00 \\ \hline
		\end{tabular}
	\end{center}
\end{table}

\subsection{ One-Dimensional Problem with Anisotropic Incoming Flux}

\begin{figure}[ht!]
\begin{center}
  \includegraphics[width=8cm]{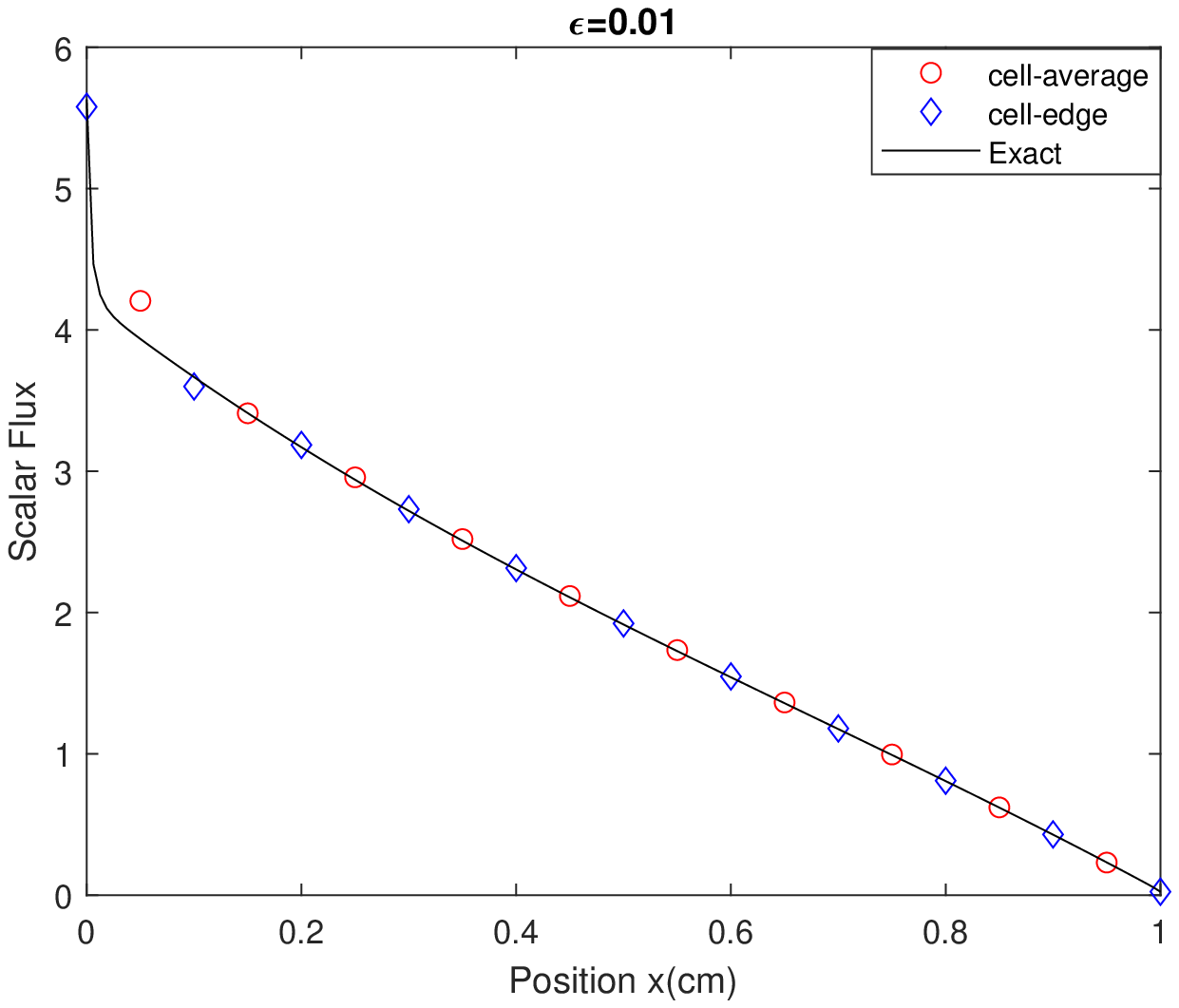}
  \includegraphics[width=8cm]{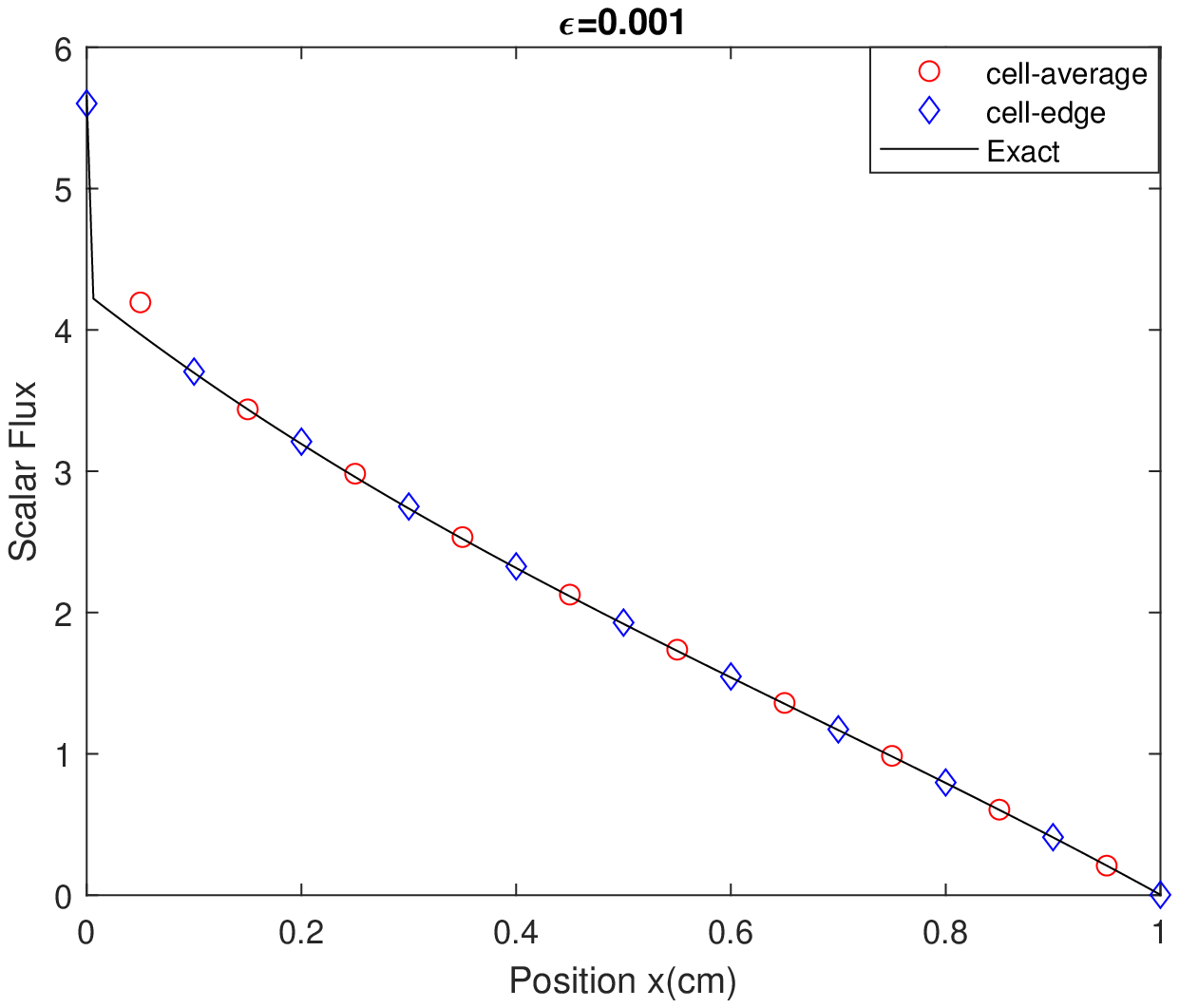}
\end{center}
 \caption{Numerical solution of Example 3 with $\Delta x=0.1$. $\varepsilon=0.01$ (left) and $\varepsilon=0.001$ (right)}\label{e3fig}
\end{figure}

\noindent\textbf{Example 3}. This is also a 1D slab case, and the Gauss-Legendre $S_{12}$
quadrature set is used for the angular discretization. The setup of the problem takes the form
$$L=1,~~~\sigma_{t}=1,~~~\sigma_{a}=0.8, ~~~ Q=1.$$
The incoming angular flux at $x=0$ changes linearly from $0$ to $5$ for the six discrete incoming directions.
On the right boundary at $x=1$, the vacuum boundary is considered.
Again, the analytical solution can be obtained following the approach discussed in \cite{wanganalytic}.
Fig. \ref{e3fig} shows the numerical solutions with spatial size $\Delta x=0.1$ for different $\varepsilon$, which demonstrates that the HWENO FSM can capture the thick
diffusion limit well in both the cell-average and cell-edge fluxes, for various values of $\eps$.

\subsection{One-Dimensional Problem with the Interior Thin Layer}
\begin{figure}
\begin{center}
  \includegraphics[width=8cm]{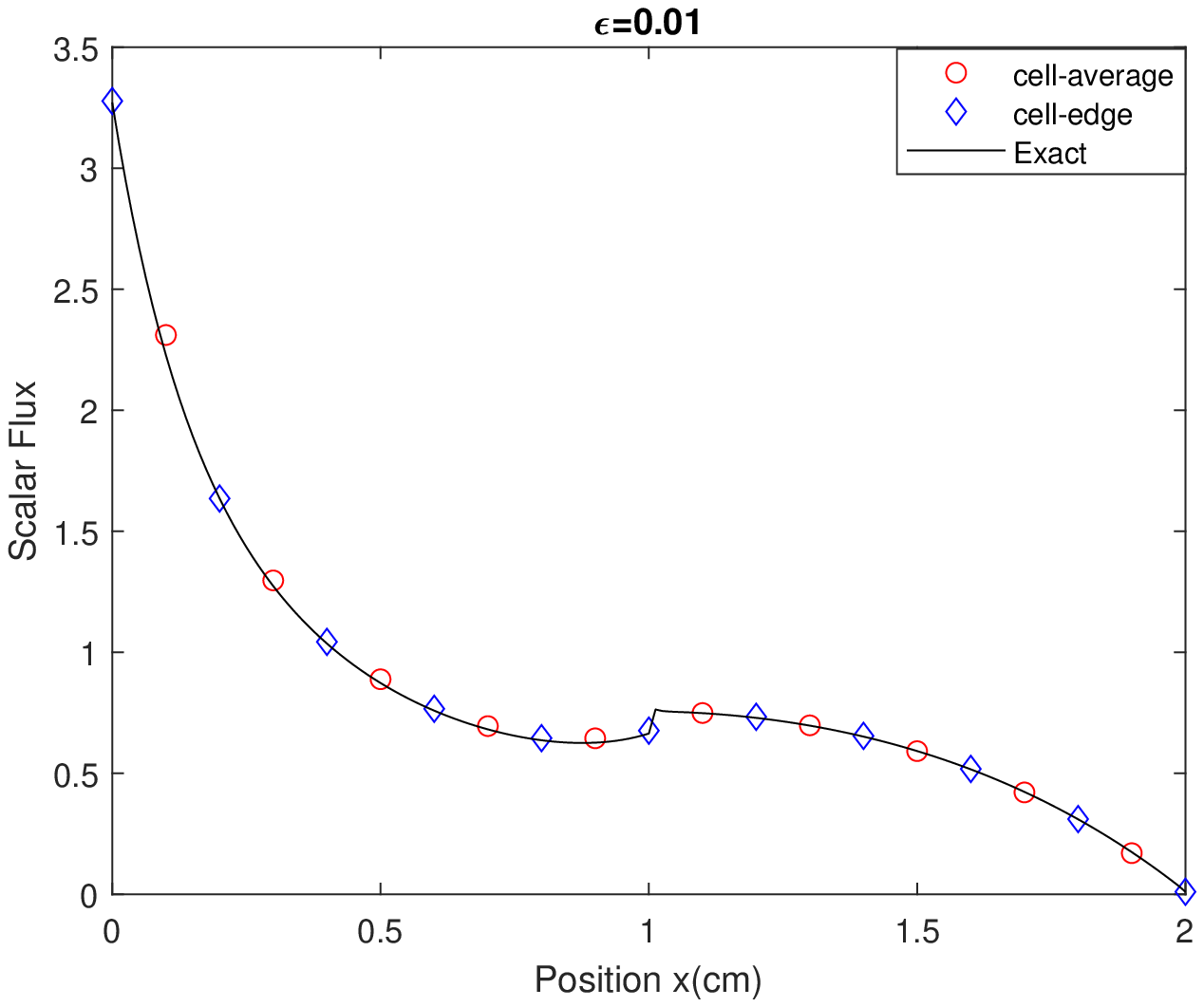}
  \includegraphics[width=8cm]{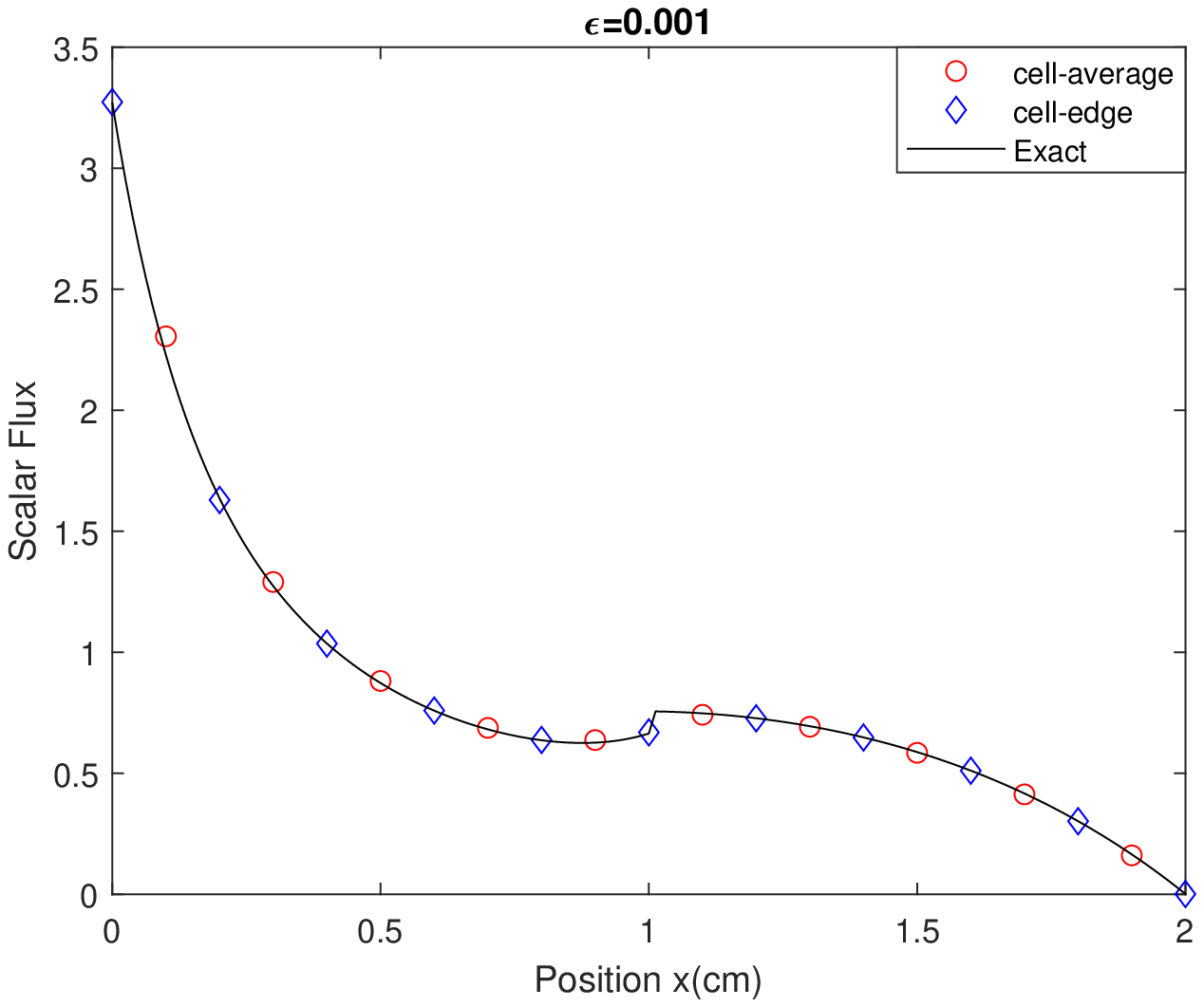}
\end{center}
 \caption{Numerical solution of Example 4 with $\Delta x=0.2$. $\varepsilon=0.01$ (left) and $\varepsilon=0.001$ (right)}\label{e4fig}
\end{figure}

\noindent \textbf{Example 4}. In this test, we consider a 1D slab consisting of two material regions. The left half of the slab is an
optically thin region, and the right half is an optically thick diffusive region.
The specifications of the problem are defined by $L=2$,
\begin{equation*}
\sigma_t=\begin{cases}\varepsilon,\quad 0\leq x<1,\\
  1,\quad 1\leq x<2,\end{cases}\quad
  \sigma_a=\begin{cases}\frac{1}{\varepsilon},\quad 0\leq x<1,\\
  0.8,\quad 1\leq x<2,\end{cases} \text{and} \quad
   Q=\begin{cases}0 ,\quad 0\leq x<1,\\
  1,\quad 1\leq x<2.\end{cases}
\end{equation*}
Again, the Gauss-Legendre $S_{12}$ quadrature set is used for the angular discretization. The incoming angular flux at $x=0$
changes linearly from $0$ to $5$ for the six discrete incoming
directions. On the right boundary at $x=2$, the vacuum boundary is considered.
The analytical solution of this problem in the form of cell-edge can be computed following the idea in \cite{wanganalytic}.
Fig. \ref{e4fig} shows the numerical solutions with spatial size $\Delta x=0.2$ for different $\varepsilon$, which demonstrates that the HWENO FSM can capture the thick
diffusion limit well in both the cell-average and cell-edge fluxes, for various values of $\eps$.

%

\begin{figure}[h!]
\begin{center}
  \includegraphics[width=8cm]{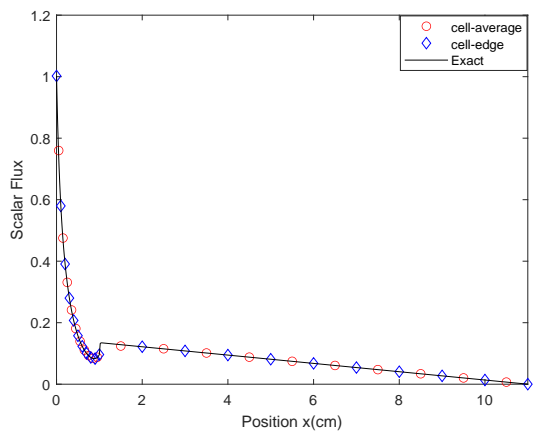}
\end{center}
 \caption{Numerical solution of Example 5 on $N=20$}\label{larsenfig1}
\end{figure}

\begin{figure}[h!]
\begin{center}
  \includegraphics[width=8cm]{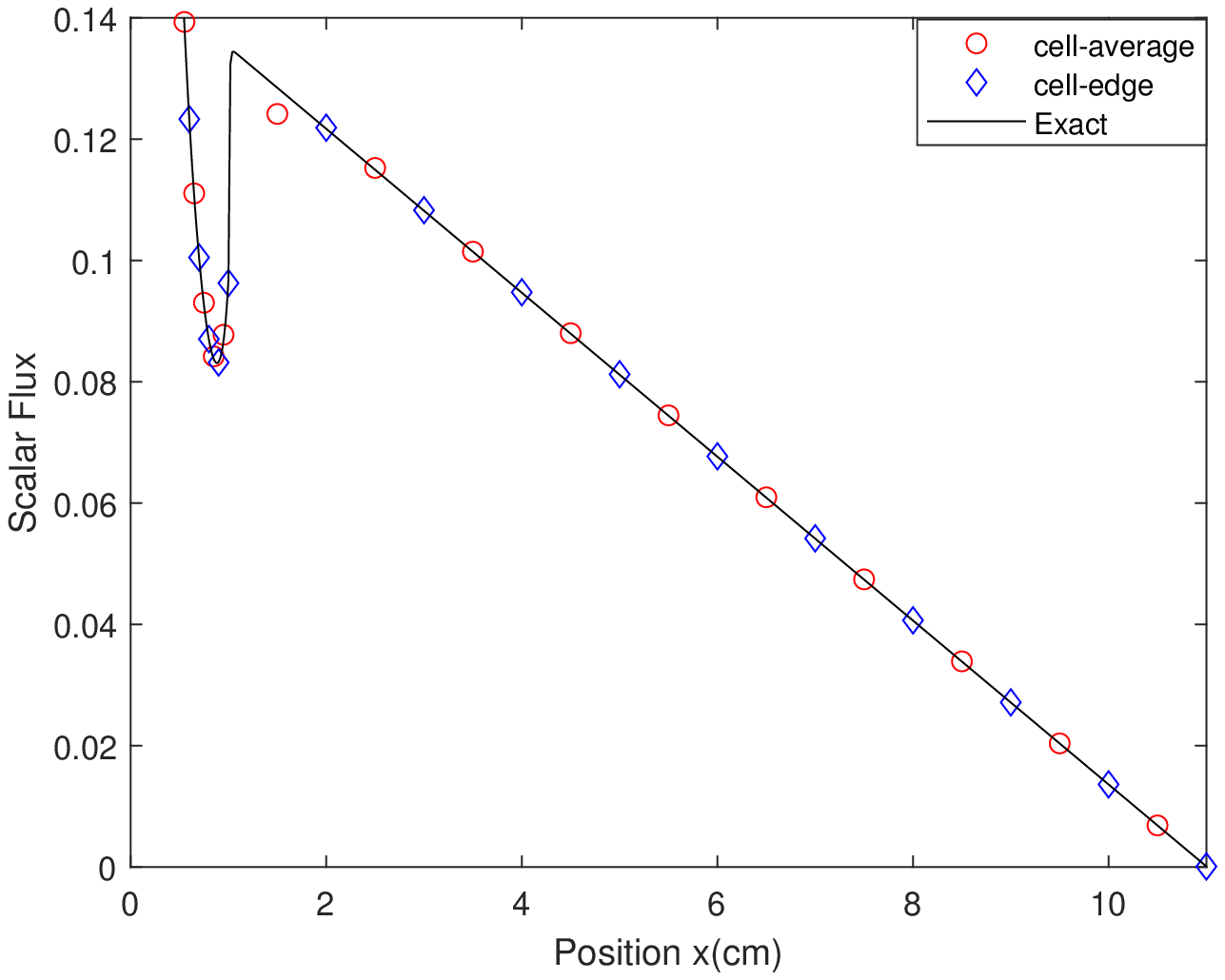}
  \includegraphics[width=8cm]{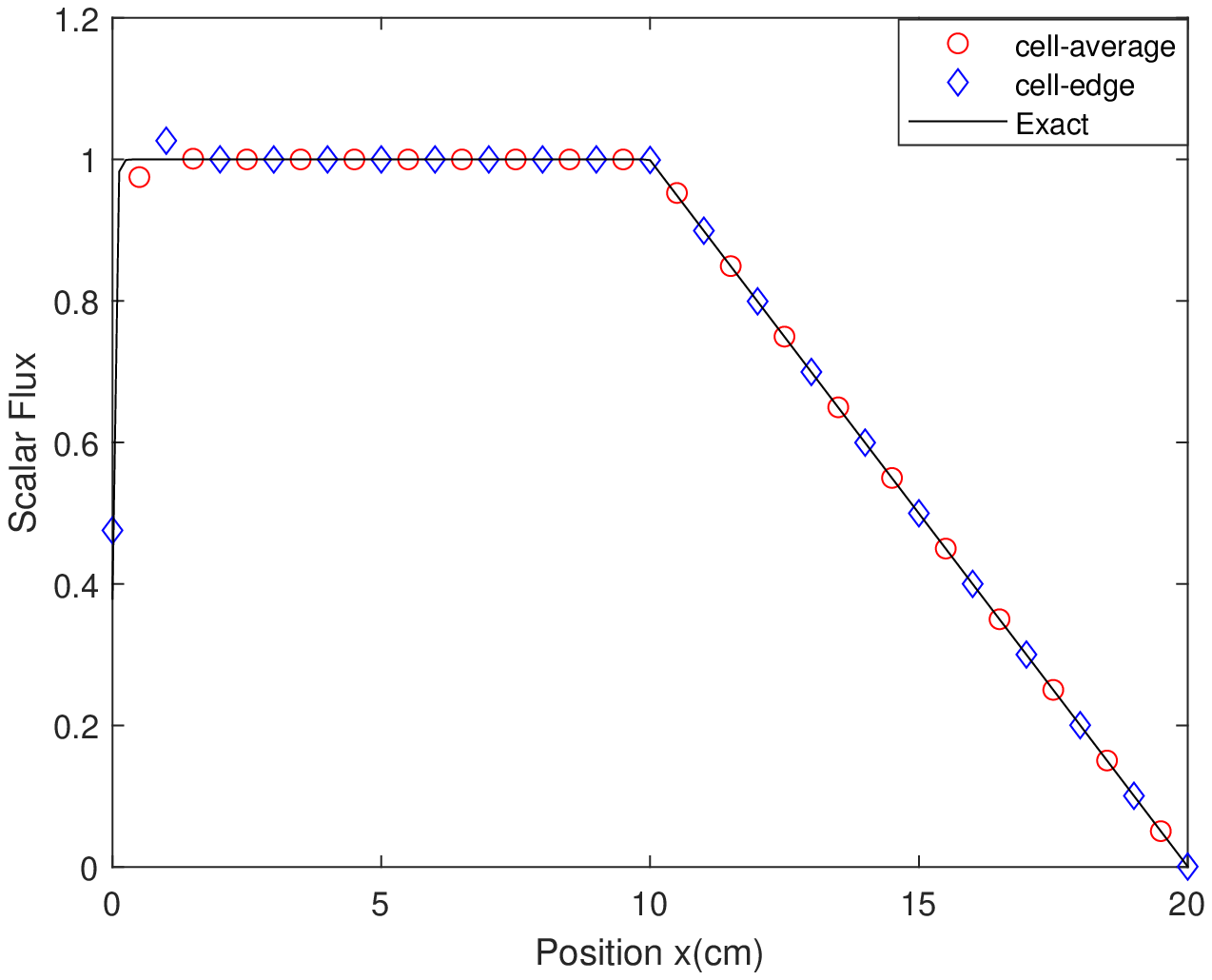}
\end{center}
 \caption{Numerical solution of scalar fluxes. left: The zoomed in version of Fig. \ref{larsenfig1}, right: Example 6 with $\Delta x=1$ }\label{larsenfig2}
\end{figure}

\noindent\textbf{Example 5}. This is a classical example taking from \cite{larsen2}.
We consider a slab in $[0,11]$, $Q=0$, $\eps=1$, and other settings are given by
 \begin{equation*}
\sigma_t=\begin{cases}2,\quad 0\leq x<1,\\
  100,\quad 1\leq x<11,\end{cases}~
  \sigma_a=\begin{cases}2,\quad 0\leq x<1,\\
  0,\quad 1\leq x<11,\end{cases} \text{and~~}
   \begin{cases}\psi(0,\mu)=1,~ \mu>0,\\
  \psi(11,\mu)=0,~\mu<0.\end{cases}
\end{equation*}
 The problem consists of a two mean-free-path purely absorbing part and a 1000 mean-free-path purely scattering part.
We solve this problem using the $S_{12}$ quadrature set in the angular direction and the mesh size
\begin{equation}
\Delta x=\begin{cases}0.1,\quad 0\leq x<1,\\
  1,\quad 1\leq x<11.\end{cases}
  \end{equation}
The ``exact" cell-edge solution is obtained by LD method with a refined mesh of $N=1000$. In Fig. \ref{larsenfig1}, we present the numerical results of HWENO method with hybrid strategy under $\varepsilon=0.01$. The zoomed in version is provided in the left plot of Fig. \ref{larsenfig2}, which provides a better view of the numerical simulation near the interior layer. We can observe that the numerical solution is in good agreement with the exact solution, which indicates that the proposed HWENO method produces very accurate result for this challenging test.

\noindent\textbf{Example 6}. Another problem considered in \cite{larsen2} has the setup
$L=20$, $\varepsilon=1$, $\sigma_{t}=100$, and the other settings are given by
\begin{equation}
  \sigma_a=\begin{cases}{10},\quad 0\leq x<10,\\
  0,\quad 10\leq x<20,\end{cases} \quad
  Q=\begin{cases}{10},\quad 0\leq x<10,\\
  0,\quad 10\leq x<20,\end{cases}
\end{equation}
with vacuum boundary. The system in this problem consists of a 1000 mean free path slab, with absorption and a flat interior source, adjoining a 1000 mean free path purely scattering slab with no interior source.
The ``exact" cell-edge solution is obtained by the LD method with a refined mesh of $N=100$. The Gauss-Legendre $S_{12}$ quadrature set is used in the angular discretization. We take spatial size $\Delta x=1$ and the numerical result of the HWENO method with $\varepsilon=0.01$ is provided in the right plot of Fig. \ref{larsenfig2}, from which we can observe that the numerical solution is in good agreement with the exact solution.

For problems 5 and 6, we also provided the computational time comparison of HWENO method with or without the hybrid strategy in Table \ref{timetable}, from which we can observe a $50\%$ saving of CPU time when the hybrid strategy is used.

\begin{table}[ht!]
\caption{The CPU times of the HWENO method with or without the hybrid strategy for all the Examples. ``iter" and ``time" denote the iteration numbers and CPU time, respectively. The ``ratio" represents the ratio of the CPU time of method with hybrid strategy over that without hybrid strategy. }\label{timetable}
	\begin{center}
		\begin{tabular}{|c|c|c|c|c|c|} 	\hline
 &\multicolumn{2}{c|}{ with hybrid strategy}	 &\multicolumn{2}{c|}{ without hybrid strategy} &		\\ \cline{1-5}		
Test& iter& time (sec)&iter&time (sec)&ratio\\ \hline
1 & - &83.95 &- &452.35 &18.55\%\\ \hline
2 & - &17.47 &- &62.73 &27.85\%\\ \hline
3 & 4726476 &123.74 &4727810 &246.88 &50.12\%\\ \hline
4 & 4456737 &153.51 &4457270 &229.46 &66.90\%\\ \hline
5 & 5311441 &209.95 &5311566 &456.40 &46.00\%\\ \hline
6 & 5680997 &221.31 &5681346 &480.77 &46.03\%\\ \hline
		\end{tabular}
	\end{center}
\end{table}

\subsection{Two-Dimensional Problems}
Only the cell-average fluxes will be plotted in the figures for the two-dimensional problems in this subsection.

\begin{figure}[ht!]
\begin{center}
  \includegraphics[width=8cm]{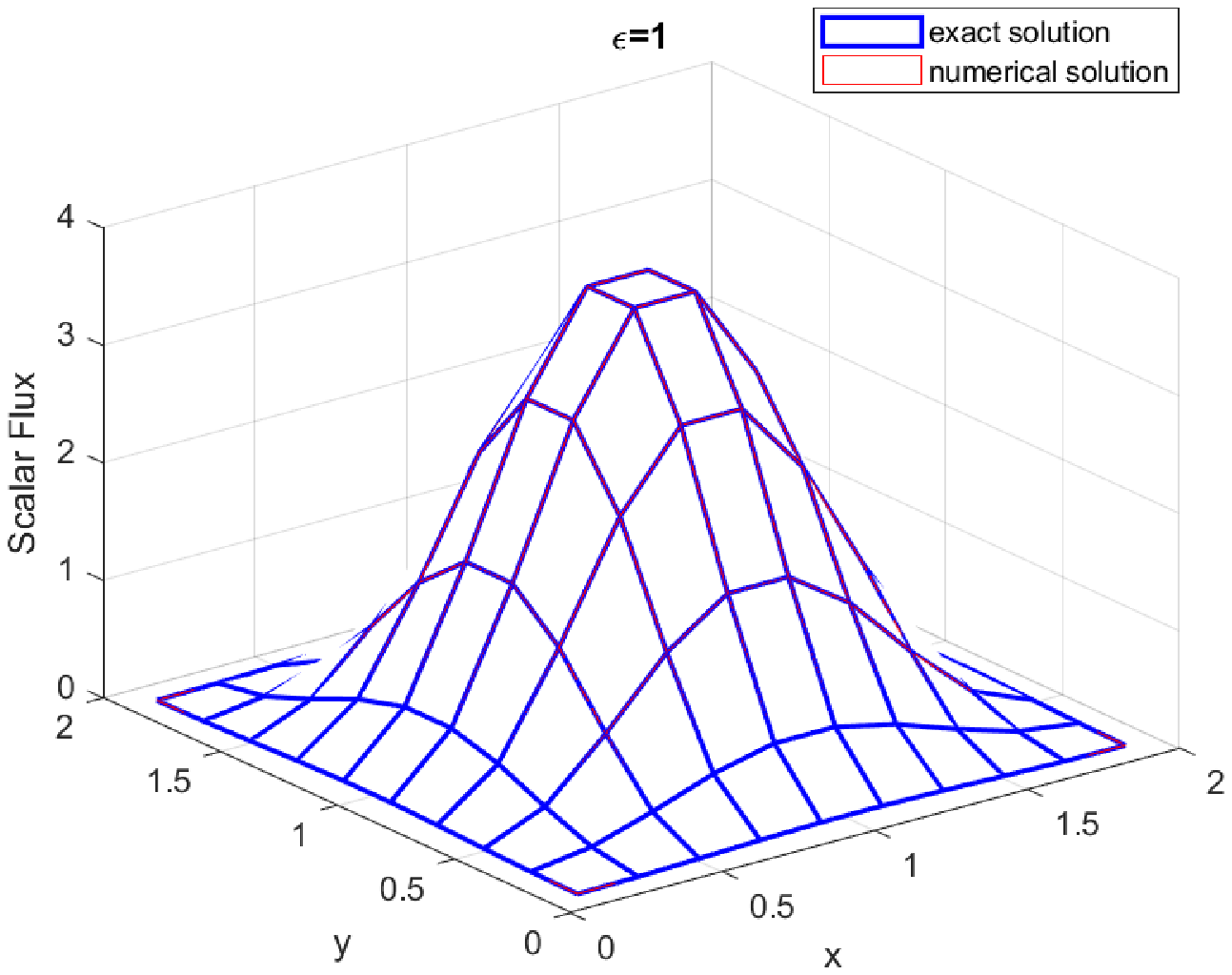}
  \includegraphics[width=8cm]{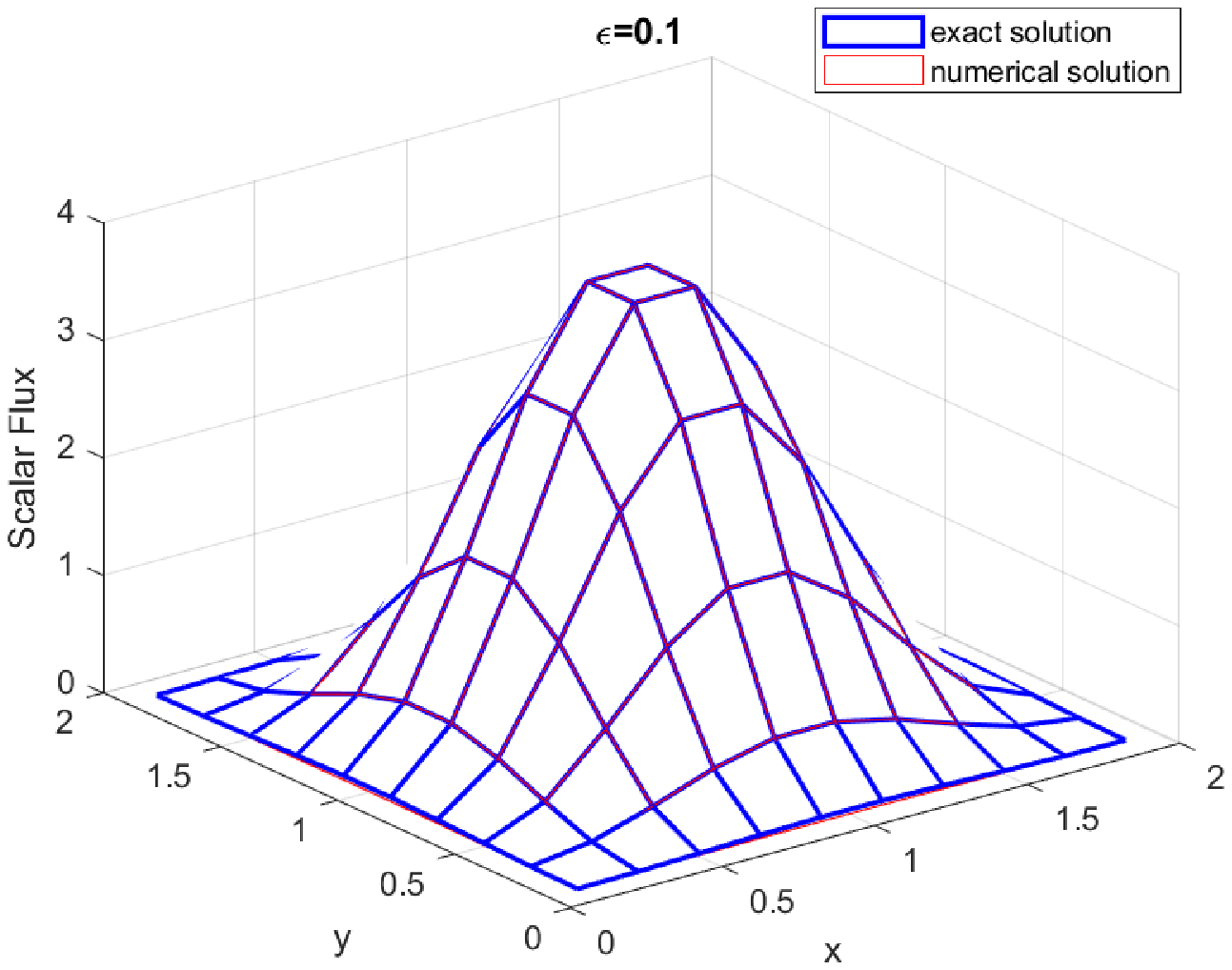}
\end{center}
 \caption{Numerical solution of Example 7 with $\Delta x=\Delta y=0.2$, and $\varepsilon=1$ (left) and $\varepsilon=0.1$ (right)}\label{e6fig1}
\end{figure}

\begin{figure}[ht!]
\begin{center}
  \includegraphics[width=8cm]{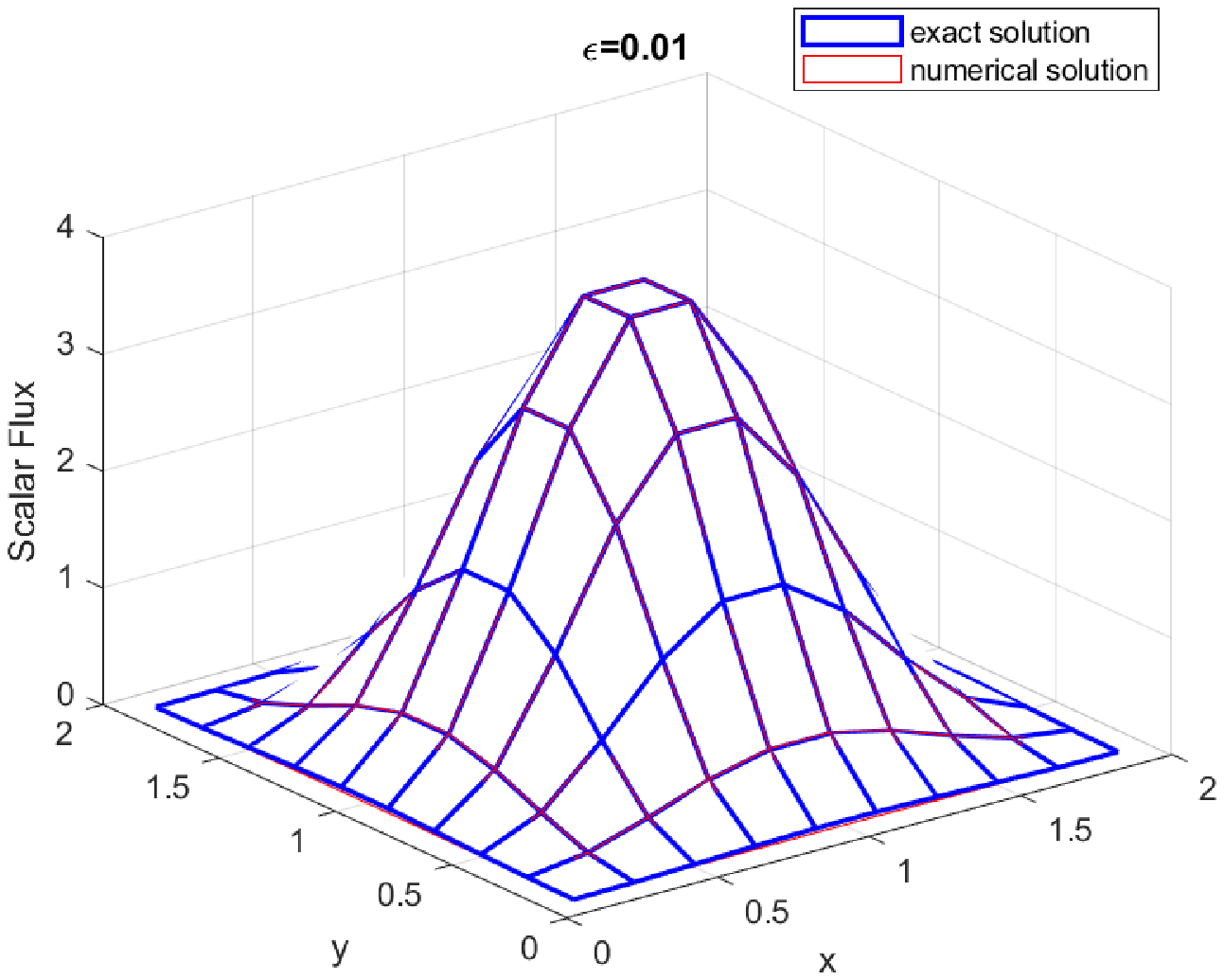}
  \includegraphics[width=8cm]{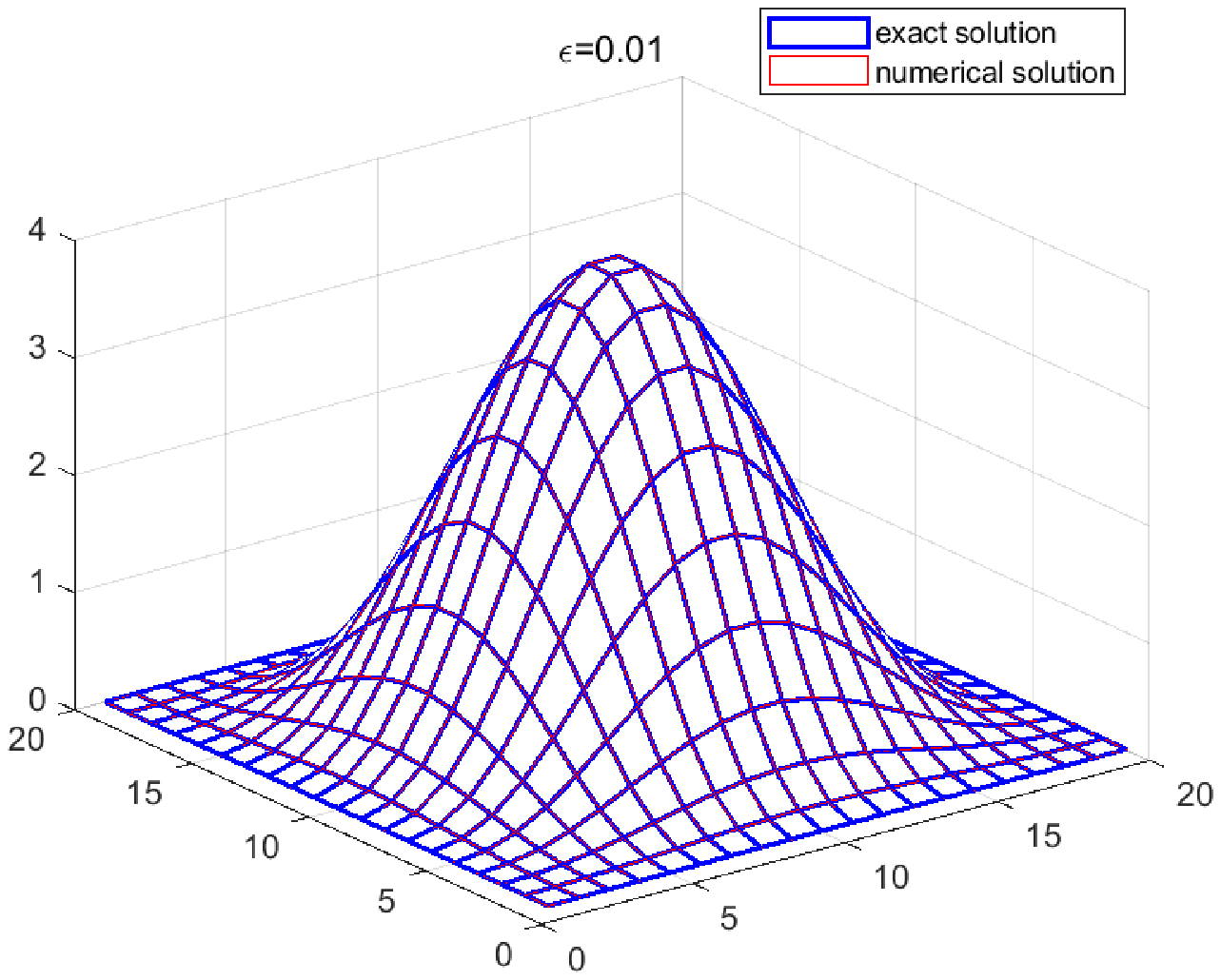}
\end{center}
 \caption{Numerical solution of Example 7 with $\varepsilon=0.01$ and $\Delta x=\Delta y=0.2$ (left) and $\Delta x=\Delta y=0.1$ (right)}\label{e6fig2}
\end{figure}

\noindent\textbf{Example 7} (Accuracy test with manufactured solution in 2D). To test the order of convergence of the 2D HWENO FSM, we follow the setup in \cite{wangap} and consider the manufactured exact solution of the form
$$\psi(x,y,\mu_{m},\eta_{m})=x^{3}y^3(2-x)^3(2-y)^3,$$
in the computational domain $\Omega=[0,2]\times [0,2]$.
The other parameters are set as
$$
\sigma_{t}=1,\qquad \sigma_{a}=0.8,$$
and
$$Q(x,y)=\frac{4}{\varepsilon}[(24x^2-48x^3+30x^4-6x^5)y^3(2-y)^3\mu_{m}+x^3(2-x)^3(24y^2-48y^3+30y^4-6y^5)\eta_{m}]+4\sigma_{a}\psi.$$
The numerical solutions are obtained using the level symmetric $S_{12}$ quadrature set for angular discretization.
We have run the simulations for $\eps=1$ and $\eps=0.1$.
In Table \ref{e6table}, we show the numerical errors, the corresponding order of accuracy and computational time of the 2D HWENO method with the hybrid strategy,
from which we can observe the expected high order accuracy for both choices of $\eps$.
In Fig. \ref{e6fig1}, we plot the numerical solutions with $\Delta x=\Delta y=0.2$ and $\eps=1$ or $\eps=0.1$.
In Fig. \ref{e6fig2}, we plot the numerical solutions with a smaller $\eps=0.01$ on two sets of computational meshes. From these figures, one could observe that the HWENO FSM can capture the thick diffusion limit well on coarse meshes.

\begin{table}
\caption{Example 7. The errors, order of accuracy and CPU time of 2D HWENO method with hybrid strategy}\label{e6table}
	\begin{center}
		\begin{tabular}{|c|c|c|c|c|c|c|c|}
			\hline
\multicolumn{7}{|c|}{$\varepsilon=1$}\\ \hline
$N$ & $L_{1}$~error&order&  $L_{\infty}$~error&order&iter &time\\\hline
10 &2.07e-03 &- &1.05e-02 &- &43 &0.47\\ \hline
20 &3.71e-05 &5.80 &3.70e-04 &4.82 &49 &2.08\\ \hline
40 &4.62e-07 &6.32 &8.73e-06 &5.40 &63 &11.04\\ \hline
80 &4.38e-09 &6.72 &1.56e-07 &5.79 &91 &73.81\\ \hline
\multicolumn{7}{|c|}{$\varepsilon=0.1$}\\ \hline
$N$ & $L_{1}$~error&order&  $L_{\infty}$~error&order&iter &time\\\hline
10 &7.93e-03 &- &3.66e-02 &- &1434 &15.92\\ \hline
20 &1.41e-04 &5.80 &1.40e-03 &4.70 &1491 &72.07\\ \hline
40 &2.04e-06 &6.11 &3.74e-05 &5.22 &1605 &313.35\\ \hline
80 &2.73e-08 &6.22 &8.43e-07 &5.473 &1831 &1475.21\\ \hline
		\end{tabular}
	\end{center}
\end{table}

\begin{figure}[ht!]
\begin{center}
  \includegraphics[width=12cm]{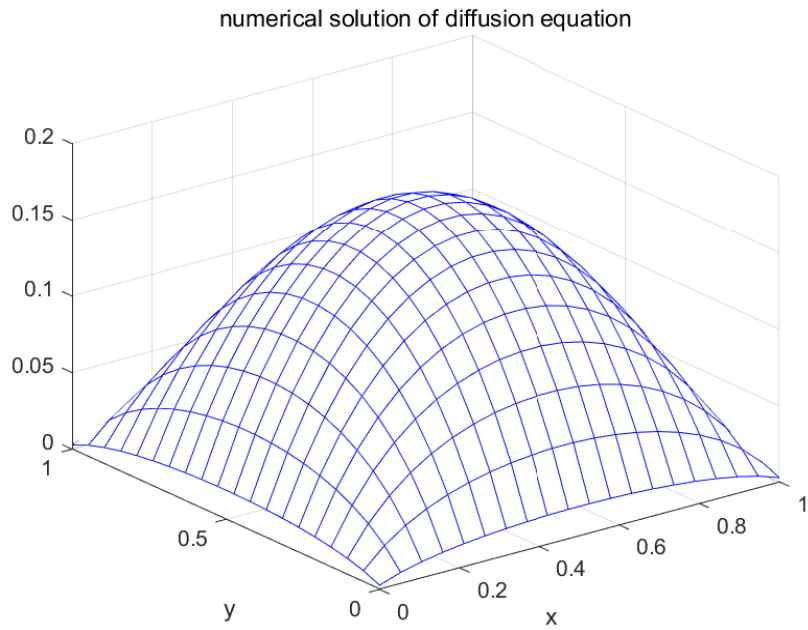}
\end{center}
\caption{Example 8: the reference exact solution of the limit diffusion equation}\label{diffusion2d}
\end{figure}

\noindent\textbf{Example 8}.  In this example which was originally considered in \cite{ld2d}, we study the problem on bounded domain with vacuum boundary conditions with the setup
$$\Omega=[0,1]^2,~
\sigma_{t}=1,~\sigma_{a}=1,~\mathrm{and}~Q=1.$$
We take the spatial size $\Delta x=\Delta y=0.05$ for this problem. Again, the level symmetric $S_{12}$ quadrature set for angular discretization. The limit diffusion equation, when $\eps \rightarrow 0$, is given by
\begin{equation}
\begin{cases}
-\frac{1}{3}\Delta\psi(\mathbf{x})+\psi(\mathbf{x})=1, &\mathbf{x}\in(0,1)^2;\\
\psi(\mathbf{x})=0, &x\in\partial(0,1)^2,
\end{cases}
\end{equation}
and we plot its exact solution on the same spatial size in Fig. \ref{diffusion2d} as a reference solution.

\begin{figure}[h!]
\begin{center}
   \includegraphics[width=8cm]{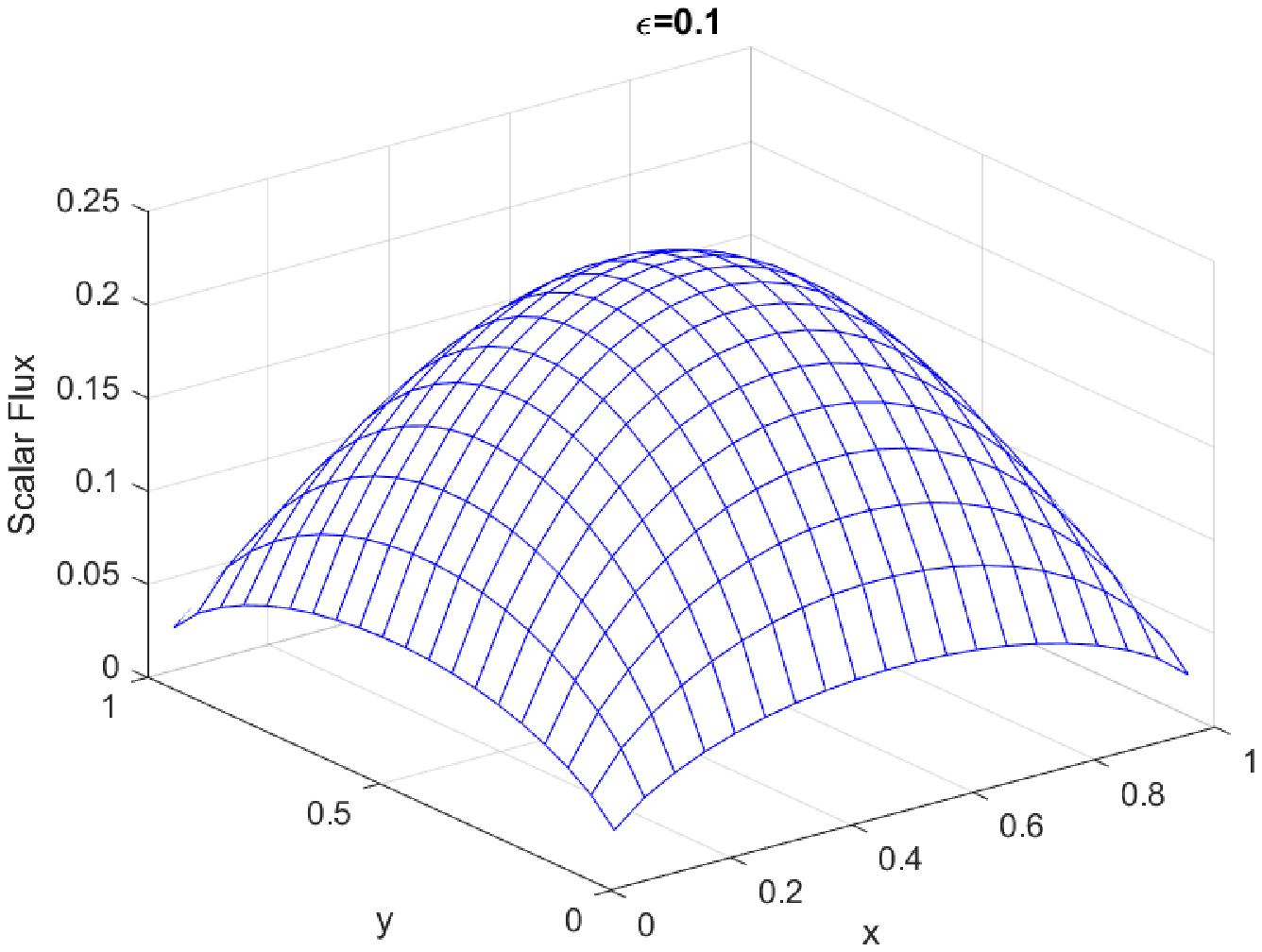}
  \includegraphics[width=8cm]{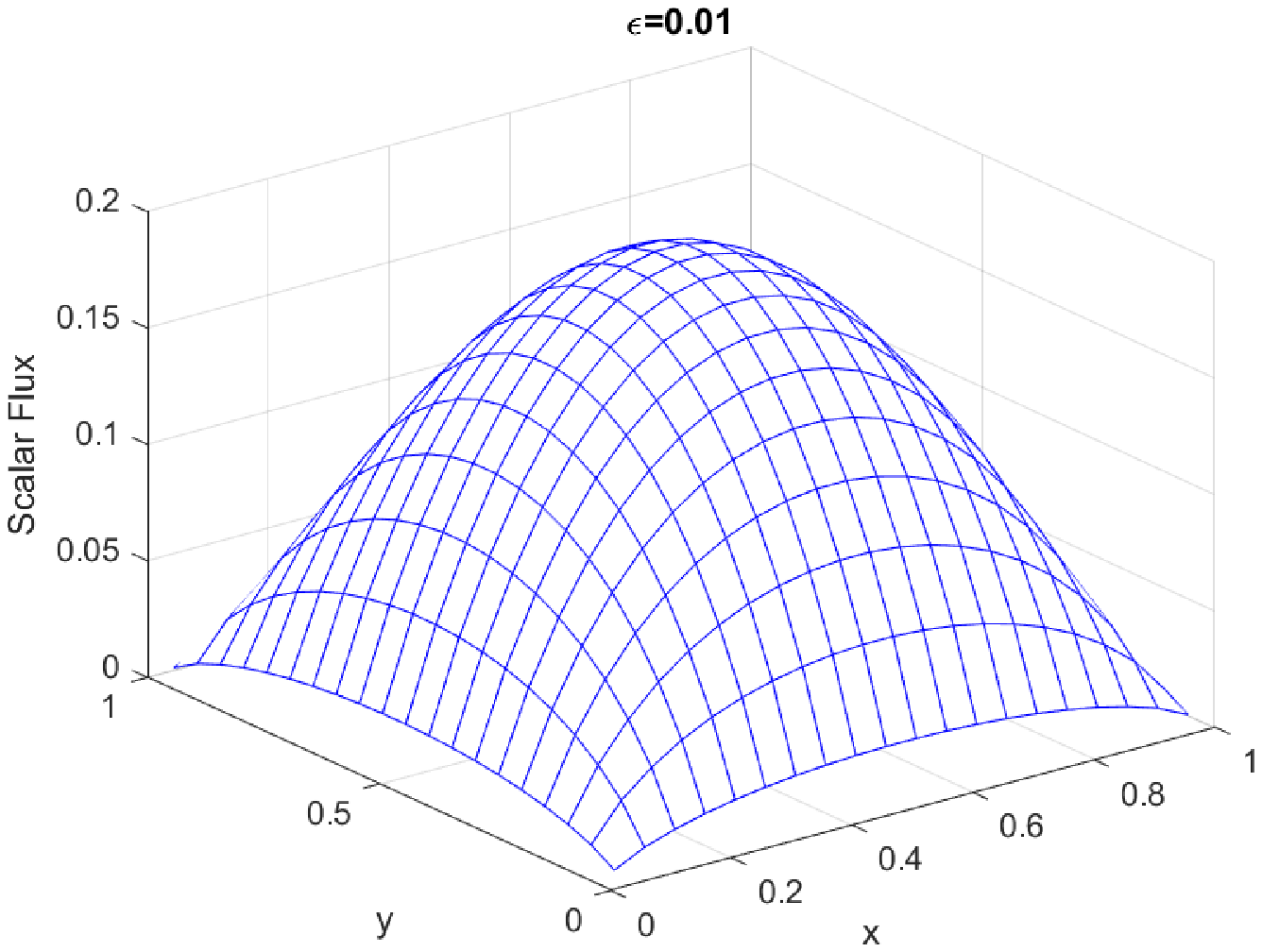}
\end{center}
 \caption{Numerical solution of Example 8 with $\varepsilon=0.1$ (left) or $\varepsilon=0.01$ (right)}\label{e8fig1}
\end{figure}

\begin{figure}[ht!]
\begin{center}
   \includegraphics[width=8cm]{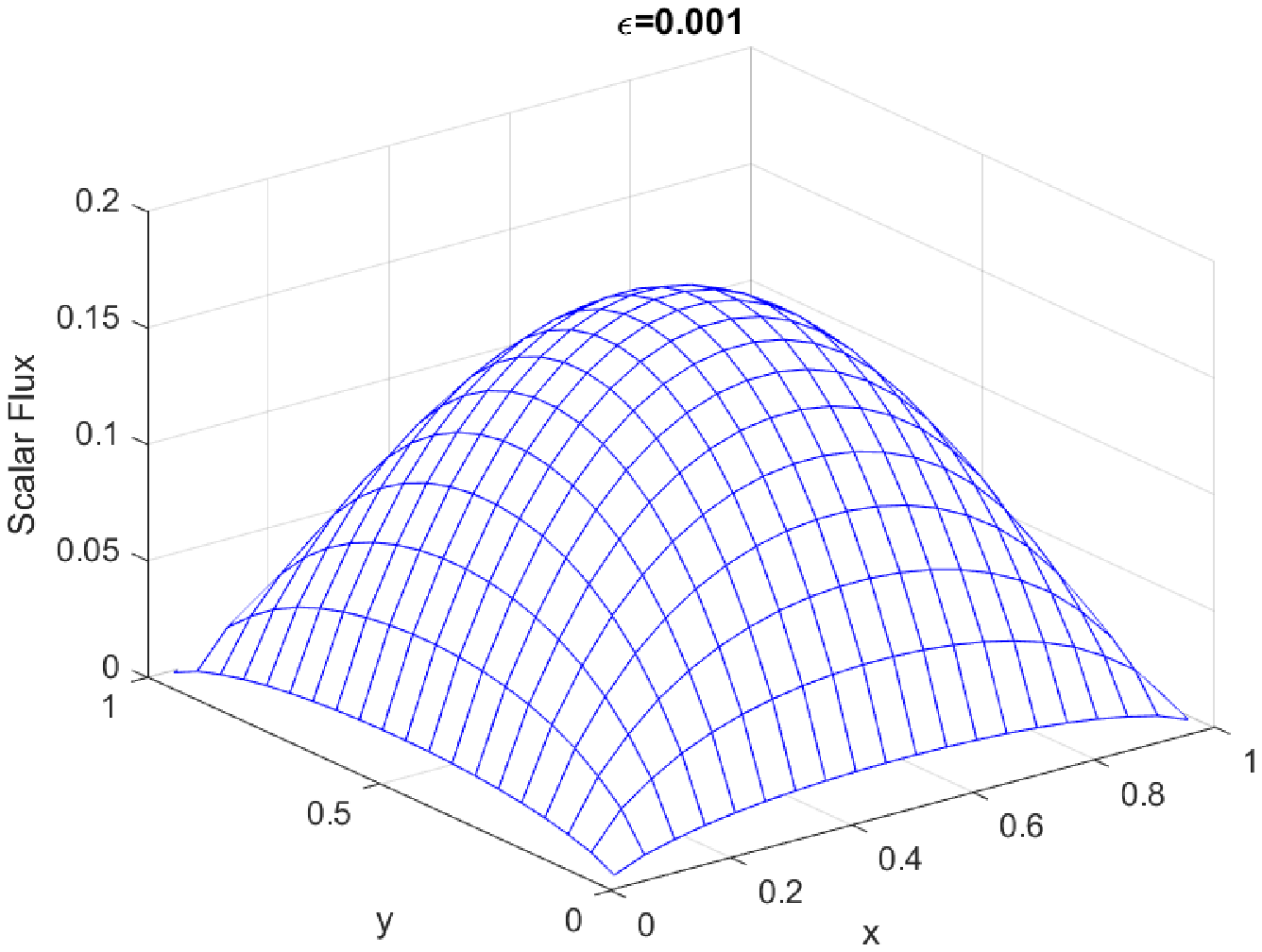}
  \includegraphics[width=8cm]{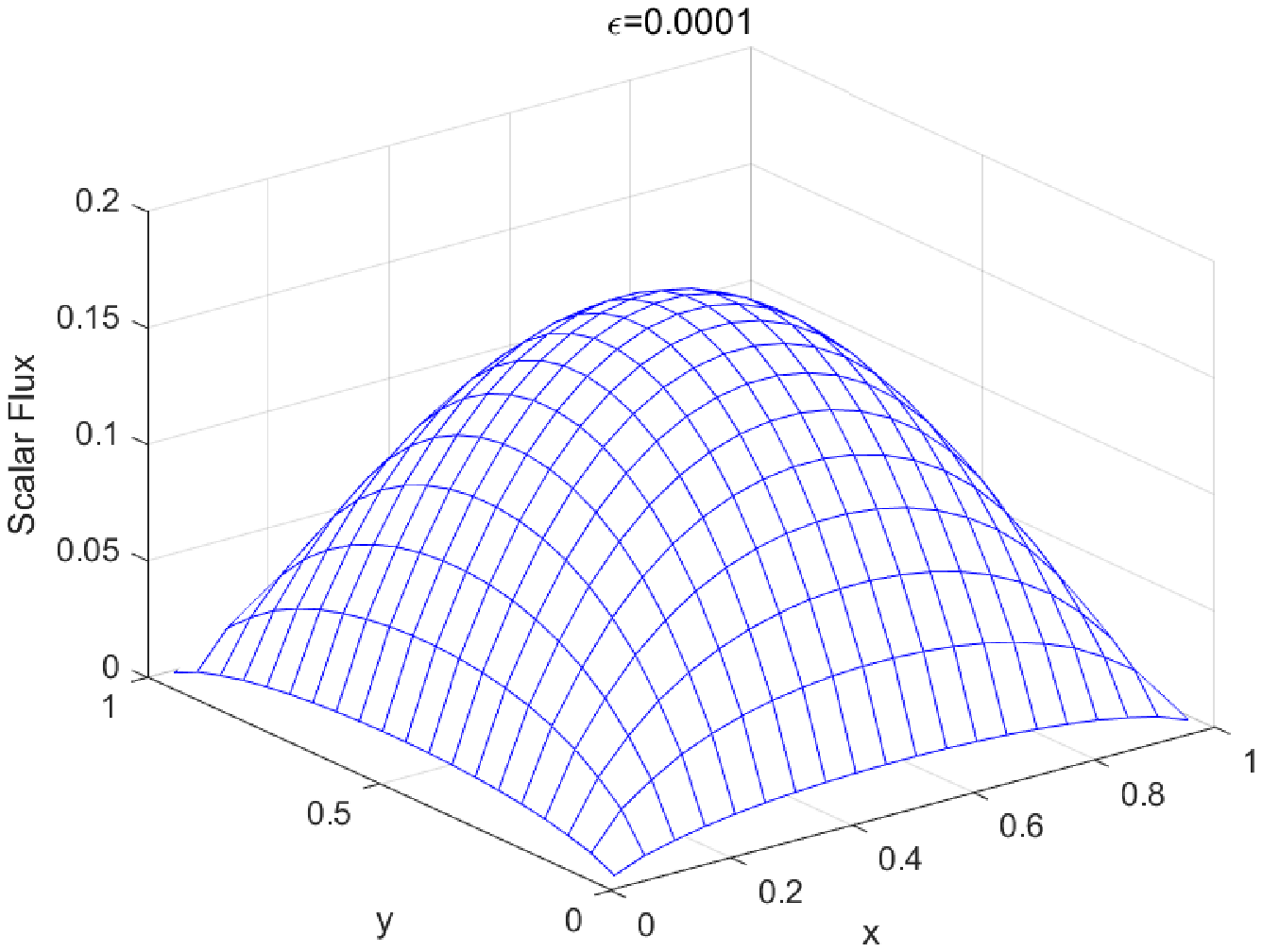}
\end{center}
 \caption{Numerical solution of Example 8 with $\varepsilon=0.001$ (left) or $\varepsilon=0.0001$ (right)}\label{e8fig2}
\end{figure}

In \cite{ld2d}, it was shown that the two-dimensional LD method does not have the AP property and cannot capture the diffusion limit well.
In Figs. \ref{e8fig1} and Fig. \ref{e8fig2}, we plot the numerical solutions of the 2D finite volume HWENO method for different values of $\eps$: $0.1$, $0.01$, $0.001$, $0.0001$, and we
can observe that the diffusion limit is well captured by our method, which suggests that HWENO method has the designed AP property.


\begin{figure}[ht!]
\begin{center}
  \includegraphics[width=12cm]{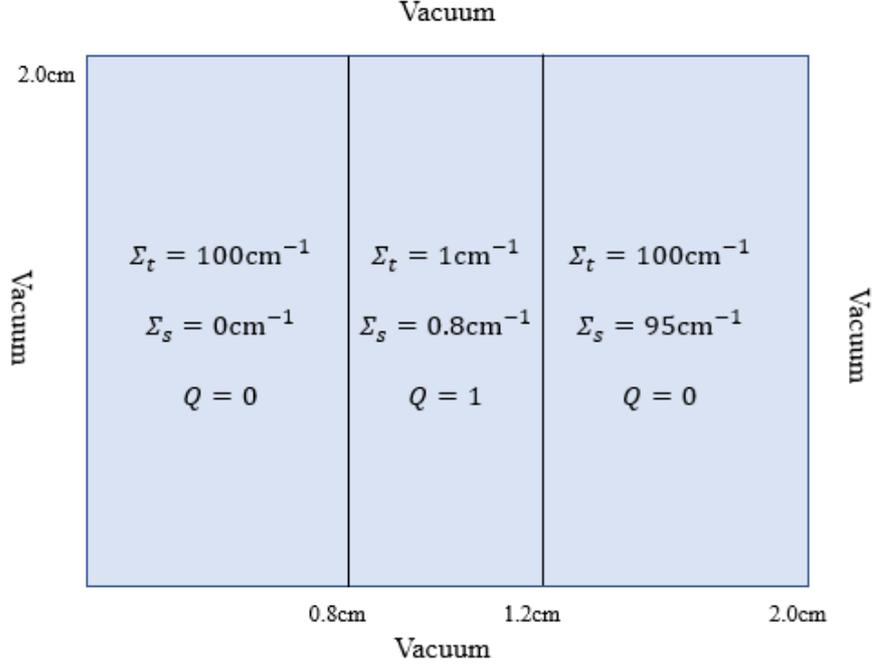}
\end{center}
\caption{Setup of Example 9}\label{e9figbc}
\end{figure}

\begin{figure}[ht!]
\begin{center}
   \includegraphics[width=8cm]{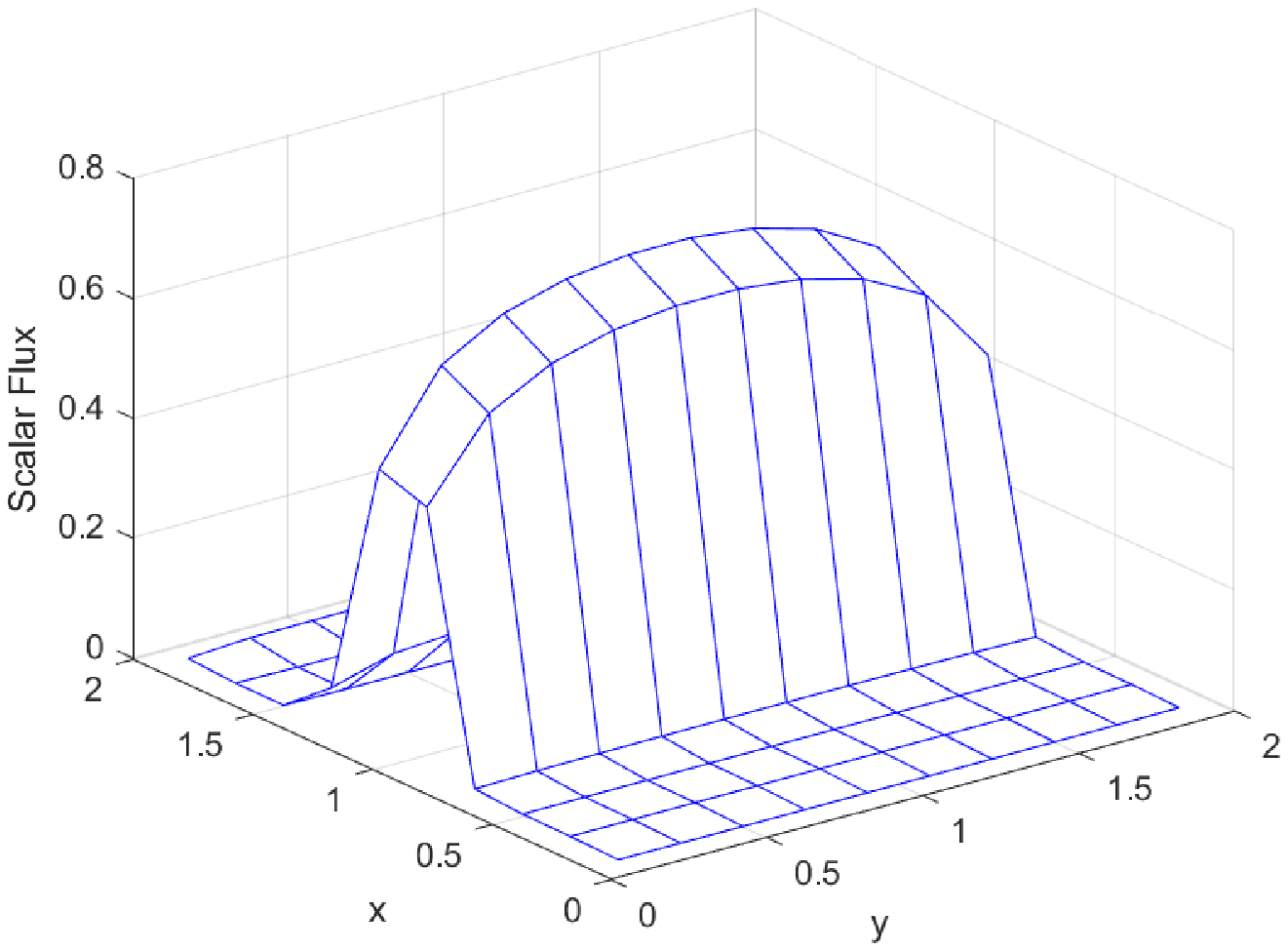}
  \includegraphics[width=8cm]{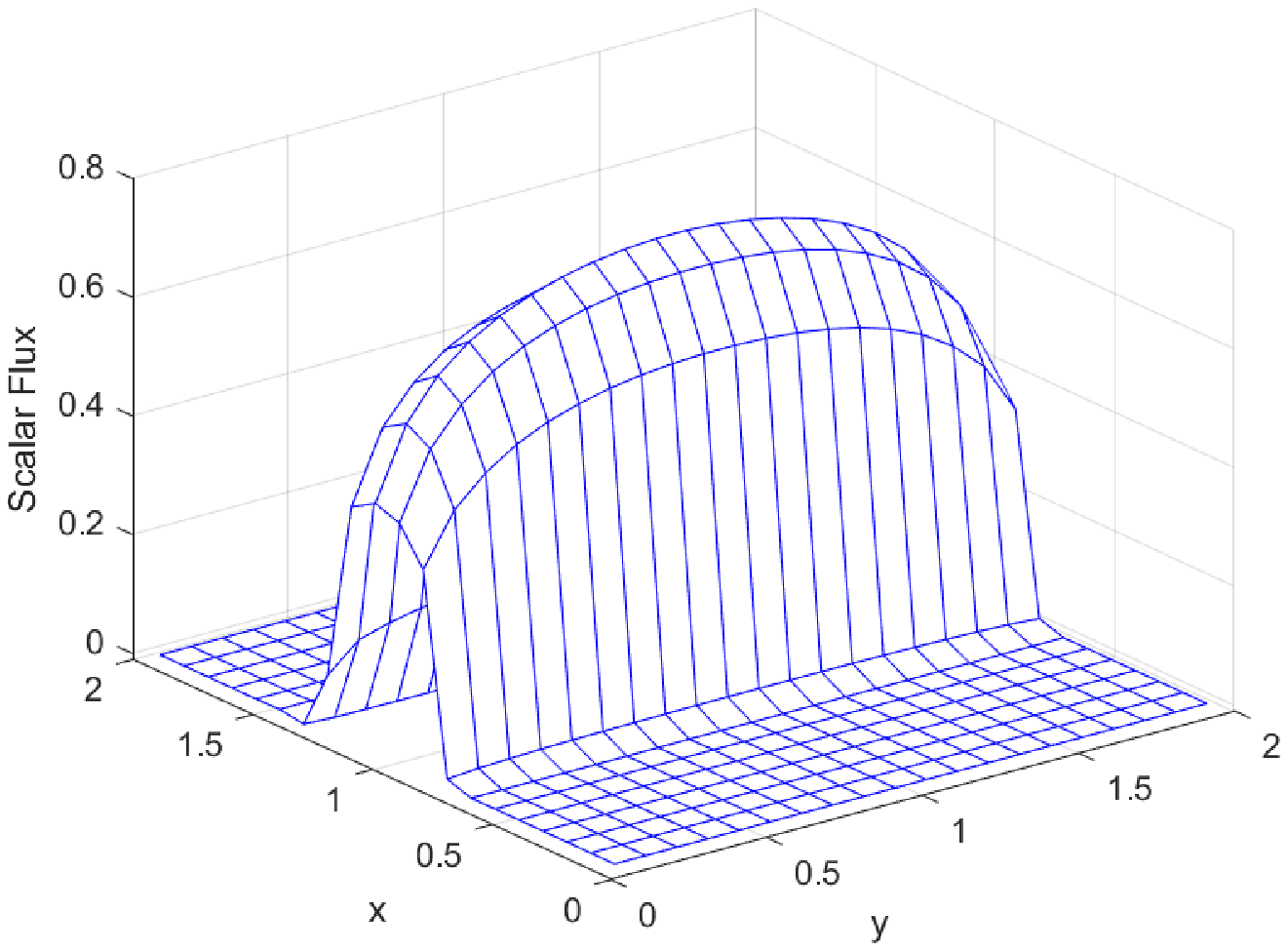}
\end{center}
 \caption{Numerical solution of Example 9 with $\Delta x=\Delta y=0.2$ (left) and $\Delta x=\Delta y=0.1$ (right)}\label{e9fig1}
\end{figure}

\noindent\textbf{Example 9}. The setup of the problem can be found in Fig. \ref{e9figbc} with $\varepsilon=1$, which consists of three subregions: the left part is a non-scattering region with no interior source, the middle part is an absorption region with interior fixed source, and the right part is a high scattering region without interior source, and a natural vacuum
boundary condition can be used for each side. The level symmetric $S_{12}$ quadrature set is used for angular discretization. The numerical results of the proposed 2D HWENO methods on different meshes are shown in Fig. \ref{e9fig1}, which again suggests that the finite volume HWENO FSM can capture the thick diffusion limit on coarse mesh.


\begin{figure}[ht!]
\begin{center}
  \includegraphics[width=8cm]{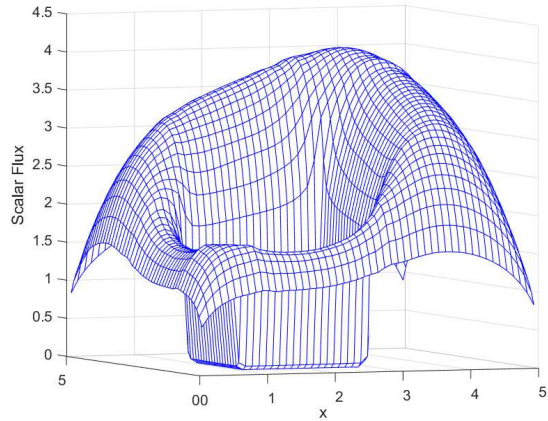}
  \end{center}
 \caption{Numerical solution of Example 10 with $\Delta x=\Delta y=0.1$. Note that the numerical solutions stay positive and there is no negative or oscillatory solution}\label{figpositive}
\end{figure}

\noindent\textbf{Example 10}.  As discussed in \cite{wanglfweno},
a potential issue with diamond difference method and other high-order $S_{N}$ numerical methods is that they could produce
non-physical negative solutions when the domain contains large material inhomogeneity.
In the last example, we consider the following problem
\begin{equation}
\Omega=[0,5]^2, \quad
\sigma_t=\begin{cases}1,\quad 0\leq x<1,\\
  100,\quad 1\leq x<3,\\
  1,\quad 3\leq x\leq5,\end{cases}\quad
  \sigma_a=\begin{cases}0.05,\quad 0\leq x<1,\\
  95,\quad 1\leq x<3,\\
  0.05,\quad 3\leq x\leq5,\end{cases}\quad ~ Q=1,
\end{equation}
 and $\varepsilon=1$. The vacuum boundary condition is considered. This problem is a 2D square problem and has a high absorbing region in the domain. The level symmetric $S_{12}$ quadrature set for angular discretization.
In Fig. \ref{figpositive}, we plot the the scalar flux distribution calculated by the proposed HWENO FSM under mesh size $\Delta x=\Delta y=0.1$,
from which we can observe the nice positivity-preserving property of the HWENO method for this example.
We refer to \cite[Fig 4]{wanglfweno} for the numerical results of diamond difference method, with negative values generated near the boundary of the
absorbing region.

\bigskip
\section{Conclusion}

In this paper, we combined the HWENO scheme and the fast sweeping method to numerically solve the steady-state $S_{N}$ transport equation in the finite volume framework. The HWENO method is known to be more compact and produce smaller error than WENO method, and enjoys a simpler boundary treatment.
The main contribution of this paper is to demonstrate that the finite volume HWENO method preserves the asymptotic limit when $\eps$ goes to zero and has the asymptotic preserving property.
One- and two-dimensional numerical examples show that HWENO FSM is of high order accurate and can capture the thick diffusion limits on coarse meshes.
The proposed method can be easily extended to any dimension on cartesian meshes. The extension to unstructured meshes, as well as comparison with other high order methods, will be discussed in future work.

\end{document}